# INTRODUCTION TO PSEUDO-DIFFERENTIAL OPERATORS

M. S. JOSHI

## 1. Preface

These notes cover most of a Part III course on pseudo-differential operators. They assume the reader is familiar with distributions particularly the Schwartz kernel theorem - the book by Friedlander provides an excellent introduction to this topic.

The point of view taken is somewhere between that of Shubin, Melrose's unpublished notes and that of Chazarain and Pirou. All of which provide good places to continue.

We assume familiarity with vector bundles and with the theory of Fredholm operators between Banach spaces. For vector bundles one source is Atiyah, [2], and for Fredholm operators the self-contained section 19.1 of Hormander Volume 3 is a good account.

## 2. Introduction

Our purpose in this course is to develop a theory of variable coefficient linear, partial differential operators on manifolds. We will then develop applications of this theory. This will allow us to "solve" elliptic operators and in particular the Laplacian acting on $k-$ forms will allow us to prove the Hodge theorem which states that the de Rham cohomology of a compact manifold is just the kernel of the Laplacian which relates the Riemannian geometry of a manifold to its topology - it will also give us an immediate proof that the de Rham cohomology of a compact manifold is finite dimensional and an easy proof of Poincare duality. A fundamental part of this will be the proof of elliptic regularity - any distributional solution of an elliptic equation is a smooth function - which will be an easy consequence of our calculus. We will also apply the calculus to the study of the propagation of singularities for non-elliptic equations - in particular we will show that associated to any linear PDE there is a Hamiltonian function on phase space that governs the propagation of singularities of solutions. This theorem gives precise information about where future singularities are, given the knowledge of where they are at some time - this presages the more advanced theories of PDEs which use symplectic geometry more and more which we will touch on, but not do in this course.

We will first study operators on $\mathbb{R}^n$ and then later shift to manifolds.

**Definition 2.1.** *A linear partial differential operator of order $m$ on $\mathbb{R}^n$ is a map from $C^\infty(\mathbb{R}^n) \to C^\infty(\mathbb{R}^n)$ of the form*

$$Pu(x) = \sum_{|\alpha| \leq m} f_\alpha(x) D_x^\alpha u(x)$$

*where $f_\alpha \in C^\infty(\mathbb{R}^n)$.*

Note

$$D_x^\alpha = \left(\frac{1}{i}\right)^{|\alpha|} \frac{\partial}{\partial x_1}^{\alpha_1} \cdots \frac{\partial}{\partial x_n}^{\alpha_n}$$

and

$$|\alpha| = \alpha_1 + \cdots + \alpha_n.$$

We put in $1/i$ as its makes dealing with Fourier transforms easier.

The fact that $\widehat{D_x^\alpha f} = \xi^\alpha \hat{f}(\xi)$ means that for $u$ of compact support, we have

$$Pu(x) = \left(\frac{1}{2\pi}\right)^n \int e^{ix.\xi} p(x,\xi)\hat{u}(\xi)d\xi, \tag{1}$$

where

**Definition 2.2.** *The total symbol of $P$ is*

$$p(x,\xi) = \sum_{|\alpha| \leq m} f_\alpha(x)\xi^\alpha = e^{-ix.\xi}P(e^{ix.\xi})$$

*and the principal symbol of $P$, denoted $\sigma_m(P)$, is*

$$p_m(x,\xi) = \sum_{|\alpha| = m} f_\alpha(x)\xi^\alpha.$$

So any differential operator can be thought of as a function on $\mathbb{R}_x^n \times \mathbb{R}_\xi^n$ - this can be thought of as phase space *or the cotangent bundle.*

To work *microlocally* means to work on phase space rather than on $\mathbb{R}^n$ and this will lead to a more coherent theory. The basic idea in micro-local analysis is that a differential operator is approximated by a multiplication operator on functions on phase space. To make sense of this idea we will introduce the algebra of pseudo-differential operators. Our starting point is (1) - $p(x,\xi)$ is a polynomial in $\xi$ and smooth in $x$ so we can think of it a finite sum of homogeneous functions $p_{m-j}(x, \xi)$ where

$$p_{m-j}(x, \lambda\xi) = \lambda^{m-j} p_{m-j}(x, \xi)$$

with $j$ running from 0 to $m$. We therefore define a (classical) pseudo-differential operator, of order $m$, to be an operator of the form (1) where $p(x, \xi)$ is an infinite sum (in a sense, we will later make precise) of terms $p_{m-j}(x, \xi)$ which are smooth in $\xi \neq 0$ and

$$p_{m-j}(x, \lambda\xi) = \lambda^{m-j} p_{m-j}(x, \xi) \ \lambda > 0.$$

We do not require smoothness up to 0 because we would be left only with polynomials! (see if you can prove this) As before $p(x, \xi)$ is the total symbol and $p_m$ the principal symbol - often denoted $\sigma_m(P)$.

This class is then denoted $\Psi DO_{cl}^m(\mathbb{R}^n)$. If we denote the space of homogeneous functions of order $m$ by $H^m$ we have that

$$\sigma_m : \frac{\Psi DO_{cl}^m}{\Psi DO_{cl}^{m-1}} \to H^m$$

is an isomorphism.

We will show that these classes form a graded algebra under composition i.e. that

$$\circ : \Psi DO_{cl}^m(\mathbb{R}^n) \times \Psi DO_{cl}^{m'}(\mathbb{R}^n) \to \Psi DO_{cl}^{m+m'}(\mathbb{R}^n).$$

The relationship of the total symbol of the product to the symbols of the original operators is complicated but at the principal level one just takes a product:

$$\sigma_{m+m'}(PQ) = \sigma_m(P)\sigma_{m'}(Q).$$

The identity operator, $I$, is a pseudo-differential operator with total symbol 1 and principal symbol 1.

So if we want to solve

$$PQ = I,$$

we need

$$1 = \sigma_0(PQ) = \sigma_m(P)\sigma_{-m}(Q).$$

So provided $\sigma_m(P)$ is never zero we take $Q_0$ to be of order $-m$ with principal symbol $\sigma_m(P)^{-1}$. We then have that

$$\sigma_0(PQ_0) = \sigma_0(I)$$

and thus that

$$PQ_0 - I = R_1 \in \Psi DO_{cl}^{-1}.$$

We now want to solve away the error term - ie we want to find $Q_1$ such that $PQ_1 = -R_1$, putting $\sigma_{-m-1}(Q_1) = -\sigma_m(P)^{-1}\sigma_{-1}(R_1)$ we get

$$P(Q_0 + Q_1) - I = R_2 \in \Psi DO_{cl}^{-2}.$$

We can now repeat this argument to get $Q_j \in \Psi DO_{cl}^{-m-j}$ such that

$$P(Q_0 + Q_1 + \cdots + Q_{N-1}) - I \in \Psi DO_{cl}^{-m-N}.$$

If we then sum (in an appropriate sense) to get $Q$, we have

$$PQ - I \in \Psi DO^{-\infty} = \bigcap_l \Psi DO^l.$$

Such a $Q$ is called a parametrix and is an inverse at the level of singularities - on a compact manifold the error will be compact and this will imply that $P$ has finite dimensional kernel and cokernel - that is the image is of finite codimension.

The above argument can be applied even if $\sigma_m(P)$ has some zeros by working away from the zeros - so one can construct an inverse away from the zero set - so the hard part of the operator is the set where $\sigma_m(P) = 0$. This is called the characteristic variety, $\mathrm{char}(P)$ and it governs much of the operators behaviour. We use these ideas to define a notion of direction and location for a singularity of a distribution $u$ - the wavefront set:

$$\mathrm{WF}(u) = \left( \bigcap_{Pu \in C^\infty} \mathrm{char}(P) \right)$$

- this is a subset of $\mathbb{R}_x^n \times \mathbb{R}_\xi^n - \{0\}$. Note that if $\phi$ is a smooth function such that $\phi u$ is smooth then we have immediately that $(x, \xi) \notin \mathrm{WF}(u)$ for $x$ with $\phi(x) \neq 0$. We later show that

$$\mathrm{singsupp}(u) = \{x : \exists \xi, \ (x, \xi) \in \mathrm{WF}(u)\}$$

so the wavefront set contains all the information that the singular support does. The theorem of Hörmander on the propagation of singularities states that if $Pu \in C^\infty$ then $\mathrm{WF}(u)$ is a union of complete integral curves of the Hamiltonian associated to the energy function $p_m$. One interpretation of this is that the operator $P$ is a quantised version of the classical system $p_m$. When $P$ is the wave operator this means that the singularities of solutions travel along geodesics in the manifold.

So our first task in this course will be to construct the class of pseudo-differential operators described above and we will then develop the applications.

## 3. Oscillatory integrals

Our objective is to make sense of integrals of the form

$$\int e^{i\phi(x,\theta)} a(x,\theta) d\theta$$

which are not absolutely convergent. We can think of these as generalisations of the Fourier transform. We will need restrictions on $\phi$, the phase, and $a$, the amplitude, for the integral to make sense.

**Definition 3.1.** $\phi \in C^\infty(\mathbb{R}_x^n \times (\mathbb{R}_\theta^k - \{0\}))$ is a phase function, if $\phi$ is homogeneous of degree one in $\theta$, is real-valued and it has no critical points, i.e.

$$d'_{x,\theta}\phi = \left( \frac{\partial\phi}{\partial x_1}, \dots, \frac{\partial\phi}{\partial x_n}, \frac{\partial\phi}{\partial\theta_1}, \dots, \frac{\partial\phi}{\partial\theta_k} \right)$$

is never zero in all coordinates simultaneously.

**Definition 3.2.** We shall say $a$ is a symbol of order $m$ if $a \in C^\infty(\mathbb{R}^n \times \mathbb{R}^k)$ and

$$|D_x^\alpha D_\theta^\beta a(x,\theta)| \leq C_{\alpha,\beta,K} <\theta>^{m-|\beta|},$$

for $x$ in $K$ compact, all $\theta$ where

$$<\theta> = (1+|\theta|^2)^{\frac{1}{2}}.$$

The class is denoted $S^m(\mathbb{R}^n; \mathbb{R}^k)$ or just $S^m$.

Note that in particular if $a(x,\theta)$ is a polynomial in $\theta$ then it is a symbol.

Note also that the estimates are really only statements about what happens at infinity. In particular, if $\chi(x,\theta)$ is compactly supported in $\theta$ then it is an element of

$$S^{-\infty} = \bigcap_m S^m.$$

**Example 3.3.** If $a(x,\theta)$ is homogeneous of degree $m$ in $\theta$ then if we smooth off near 0, we obtain a symbol of order $m$.

The following proposition is trivial but important.

**Proposition 3.4.** Pointwise multiplication defines a map

$$S^m \times S^{m'} \to S^{m+m'}$$

and

$$D_x^\alpha D_\theta^\beta : S^m \to S^{m-|\beta|}.$$

So now lets return to our integral. For simplicity, we will assume that $a$ is zero in a neighbourhood of $\theta = 0$ but this is not an important issue as we shall see that important effects come from infinity. If $a \in S^m$ and $m < -k$ then by the dominated convergence theorem we have that the integral converges to a function which is continuous in $x$. If we differentiate the integrand with respect to $x$, we obtain an integral of the same form except that the amplitude is in $S^{m+1}$. Hence we deduce that if $m + l < -k$ then the integral yields a function in $C^l$. So if $a \in S^{-\infty}$ then the integral is smooth.

What do we do in general? We regularize the integral to obtain a distribution.

We want to make sense of the expression:

$$< u, \psi > = \iint e^{i\phi(x,\theta)} a(x, \theta) \psi(x) d\theta dx.$$

for any $\psi \in C_0^\infty(\mathbb{R}^n)$. As this is important we shall do so in two ways.

First consider the operator, $L$ defined by

$$L = \sum_{j=1}^n \frac{\partial \phi}{\partial x_j} \frac{\partial}{\partial x_j} + \sum_{l=1}^k \frac{\partial \phi}{\partial \theta_l} |\theta|^2 \frac{\partial}{\partial \theta_l}.$$

$L$ has the property that

$$L e^{i\phi} = \left( \sum_{j=1}^n \left( \frac{\partial \phi}{\partial x_j} \right)^2 + \sum_{l=1}^k |\theta|^2 \left( \frac{\partial \phi}{\partial \theta_l} \right)^2 \right) e^{i\phi}.$$

Note that our hypothesis guarantees that the coefficient will not vanish off $\theta = 0$ and we have also that is homogeneous of degree 2 from the homogeneity of $\phi$. Call the coefficient $\chi$ then putting $M = \chi^{-1} L$. We have, at least formally,

$$\iint e^{i\phi(x,\theta)} a(x, \theta) \psi(x) d\theta dx = \iint M^r e^{i\phi(x,\theta)} a(x, \theta) \psi(x) d\theta dx.$$

Formally, we can integrate by parts to obtain

$$\iint e^{i\phi(x,\theta)} (M^t)^r (a(x, \theta) \psi(x)) d\theta dx \qquad (2)$$

and

$$M^t = \sum_{j=1}^n b_j \frac{\partial}{\partial x_j} + \sum_{l=1}^k c_l \frac{\partial}{\partial \theta_l} + d$$

with $d, b_j \in S^{-1}$ and $c_l \in S^0$. (They may be singular at zero but our integrand is supported away from zero so this is irrelevant. )

The advantage of all this is that

$$(M^t)^r(a(x,\theta)\psi(x)) = \sum_{|\alpha| \leq r} a_\alpha(x,\theta)D_x^\alpha \psi(x)$$

with $a_\alpha \in S^{m-r}$ and so if $m - r < -k$ , the integral (2) will converge absolutely. So we take (2) to be the definition of

$$\int e^{i\phi(x,\theta)}a(x,\theta)d\theta$$

and it is clear that we have a distribution of order $r$. Note that our formal computation is valid if $m - r < -k$ and so the answer will be independent of the choice of $r$.

Since this is so important to everything we do and it may appear a little artificial, let's consider another approach as well. Let $\chi$ be a smooth bump function i.e.

$$\chi(\theta) = \begin{cases} 1 \text{ for } |\theta| < 1 \\ 0 \text{ for } |\theta| > 2. \end{cases}$$

Then we define

$$< u, \psi > = \lim_{\epsilon \to 0} \iint e^{i\phi(x,\theta)}\chi(\epsilon\theta)a(x,\theta)\psi(x)d\theta dx.$$

The integral is compactly supported for any fixed $\epsilon > 0$ so there is no problems with it. We want to see what happens as $\epsilon$ goes to zero. If $a \in S^m$ and $m < -k$ then the integral is uniformly absolutely convergent and so will converge to the convergent integral

$$\iint e^{i\phi(x,\theta)}a(x,\theta)\psi(x)d\theta dx.$$

So as before the issue is what happens for $m \geq -k$. Letting $M$ be as above, we have

$$< u, \psi > = \lim_{\epsilon \to 0} \iint e^{i\phi(x,\theta)}(M^t)^k \left(\chi(\epsilon\theta)a(x,\theta)\psi(x)\right) d\theta dx.$$

Without loss of generality we can assume $\epsilon$ is within a compact set so we have that

$$|D_\theta^\alpha \chi(\epsilon\theta)| \leq C_\alpha < \theta >^{-|\alpha|}$$

independently of $\epsilon$. The important point is that this estimate is uniform in $\epsilon$ all the way to $\epsilon = 0$ - for fixed $\epsilon \neq 0$ the function $\chi(\epsilon\theta)$ is, in fact, in $S^{-\infty}$ but not uniformly. Thus as before we can deduce

$$(M^t)^r \left(\chi(\epsilon\theta)a(x,\theta)\psi(x)\right) = \sum_{|\alpha| \leq k} a_\alpha(\epsilon, x, \theta)\partial_\alpha\psi$$

with $a_\alpha \in S^{m-r}(\mathbb{R}_\epsilon \times \mathbb{R}_x^n; \mathbb{R}^k)$ i.e. the behaviour of $a_\alpha$ at $\infty$ is uniform in $\epsilon$. And so taking $r$ such that $m - r < -k$ we see that the integral is uniformly convergent in $\epsilon$ and must converge to

$$\iint e^{i\phi(x,\theta)}(M^t)^r a(x,\theta)\psi(x)d\theta dx,$$

as the derivatives falling on $\chi$ will disappear when $\epsilon \to 0$. So the two approaches yield the same answer! Note for the functional analysts, one can view this as the extension of a Frechet-space-valued function from a dense subset.

In general, we will deal fairly freely with oscillatory integrals but it is important to realize that when we see such an expression that this is what we mean.

As oscillatory integrals yield distributions rather than smooth functions, one thing we want to do is find their singularities.

**Proposition 3.5.** *If*

$$u = \int e^{i\phi(x,\theta)}a(x,\theta)d\theta$$

*then*

$$\text{singsupp}(u) \subset \left\{x : \exists \theta \ (x,\theta) \in supp(a) \ and \ d'_\theta\phi(x,\theta) = 0\right\}$$

Heuristically, this says that if the phase is oscillating near $x$ then the behaviour for large $|\theta|$ cancels and so there is no singularity. This is sometimes called *stationary phase* - or the WKB method.

*Proof.* Suppose

$$d'_\theta\phi(x_0,\theta) \neq 0, \forall\theta$$

then if we multiply by a bump function, $\psi$, supported sufficiently close to $x_0$ we have

$$\psi u = \int e^{i\phi(x,\theta)}b(x,\theta)d\theta$$

with $b$ a symbol such that $d'_\theta\phi$ is non-zero on the support of $b$.

This means, similarly to above, we can put

$$M = |d_\theta' \phi|^{-2} \sum_j \frac{\partial \phi}{\partial \theta_j} \frac{\partial}{\partial \theta_j}$$

and $M$ will be non-singular on the support of $b$ and will, as before, have the property

$$M e^{i\phi} = e^{i\phi}.$$

Thus if $f \in C_0^\infty$ we have

$$< \psi u, f > = \lim_{\epsilon \to 0} \int e^{i\phi(x,\theta)} (M^t)^r (b(x,\theta)\chi(\epsilon\theta)) f(x) d\theta dx.$$

The important point here is that the derivatives do not fall on $f$ as it is independent of $\theta$. Any terms falling on $\chi$, will have coefficients which are powers of $\epsilon$ and so will disappear in the limit. Thus,

$$\psi u = \int e^{i\phi(x,\theta)} (M^t)^r b(x,\theta) d\theta.$$

But $(M^t)^r b(x,\theta) \in S^{m-r}$ so we conclude that $\psi u$ is smooth as this is true for any $r$. $\qquad \square$

**Example 3.6.** Let $\phi(x,y,\xi) = \langle x - y, \xi \rangle$ then $\phi$ is a phase function as $d_x'\phi = \xi$ is non-zero for any non-zero $\xi$ and

$$d_\xi' \phi = x - y$$

so the associated oscillatory integrals are singular only on the diagonal.

## 4. Pseudo-differential Operators

**Definition 4.1.** *We define $P$ to be a pseudo-differential operator on $\mathbb{R}^n$ of order $m$, if it has a Schwartz kernel*

$$K(x,y) = \int e^{i\langle x-y, \xi \rangle} a(x,y,\xi) d\xi$$

*with $a \in S^m(\mathbb{R}_x^n \times \mathbb{R}_y^n; \mathbb{R}_\xi^n)$. We denote this class $\Psi DO^m(\mathbb{R}^n)$.*

This is slightly different from what we defined before but we shall see that this definition is equivalent.

We want our operators to act on distributions so first we must show that they preserve the class of smooth functions.

**Theorem 4.2.** *If $P$ is a pseudo-differential operator then*

$$P : C_0^\infty(\mathbb{R}^n) \to C^\infty(\mathbb{R}^n)$$

*and is continuous and if $P$ is properly supported*

$$P : C_0^\infty(\mathbb{R}^n) \to C_0^\infty(\mathbb{R}^n).$$

*Proof.* Let $u$ be a smooth function then

$$Pu = \iint e^{i\langle x-y,\xi\rangle} a(x,y,\xi) u(y) d\xi dy.$$

For fixed $x$, this is a well-defined oscillatory integral in $(y,\xi)$ and so, since the estimates are locally uniform in $x$, this gives a continuous function of $x$. The continuity is immediate from the fact that the oscillatory integral will be estimated by a finite number of derivatives in $y$ of $u$.

The formal $x$ derivatives are of the same form and so will also be continuous.

To check that the formal derivatives and the actual ones agree we evaluate the quotient - wlog we consider differentiation in $x_1$

$$\frac{Pu(x+he_1) - Pu(x)}{h} =$$

$$\iint \frac{e^{i<x+he_1,\xi>}a(x+he_j,y,\xi) - e^{i<x,\xi>}a(x,y,\xi)}{h} e^{-i<y,\xi>} u(y) d\xi dy$$

$$=$$

$$\iint \frac{\partial}{\partial s_1} \left( e^{i(<s_1,\xi_1> + <x'',\xi''>)} a(s_1,x'',y,\xi) \right)_{|s=x_1+\epsilon h} e^{-i<y,\xi>} u(y) d\xi dy$$

with $0 \le \epsilon \le 1$. This is just the value of the formal derivative at $x_1 + \epsilon h$ and as $h \to 0$ this must converge to the value of the formal derivative at $x$ and so the result follows. (just use induction for the higher terms). $\qquad\square$

To define the action of pseudo-differential operators on distributions we need to find their adjoints. However if $K(x,y)$ is the Schwartz kernel of $P$ then $P^t$ has the kernel $K(y,x)$. So $P^t$ has a kernel

$$\int e^{i<y-x,\xi>} a(y,x,\xi) d\xi$$

and performing the change of variables $\eta = -\xi$ we have

$$\int e^{i<x-y,\eta>} a(x,y,-\eta) d\eta.$$

Thus $P^t$ is a pseudo-differential operator and is continuous on $C_0^\infty(\mathbb{R}^n)$. Thus we have

**Theorem 4.3.** *If $P \in \Psi DO^m(\mathbb{R}^n)$ and is properly supported then $P$ induces a map on $\mathcal{D}'(\mathbb{R}^n)$ via*

$$< Pu, \phi > = < u, P^t\phi >, u \in \mathcal{D}'(\mathbb{R}^n), \phi \in C_0^\infty(\mathbb{R}^n).$$

We have that properly supported pseudo-differential operators map distributions and preserve the space of compactly supported smooth functions. It is therefore not surprising that they preserve the class of smooth functions (which is a subclass of distributions.): if $P$ is properly supported and $u$ is smooth then $Pu$ is a distribution. Now given any point $x_0$ we can put $u = u_1 + u_2$ with $u_1 \in C_0^\infty$ and $u_2$ such that $x_0 \notin \text{supp}(P) \circ \text{supp}(u_2)$ which implies $x_0 \notin \text{supp}(Pu_2)$. We thus have that $Pu$ equals $Pu_1$ near $x_0$ and thus is smooth near $x_0$.

Since we have shown that pseudo-differential operators preserve smooth functions and are singular only the diagonal, it is now easy to show they have the additional property of pseudo-locality.

**Definition 4.4.** *The operator $P$ is pseudo-local if*

$$\text{singsupp}(Pu) \subset \text{singsupp}(u)$$

*for any $u \in \mathcal{D}'(\mathbb{R}^n)$.*

**Theorem 4.5.** *Properly supported pseudo-differential operators are pseudo-local.*

*Proof.* Let $P$ be a pseudo-differential operator with Schwartz kernel $K$ and $u \in \mathcal{D}'(\mathbb{R}^n)$.

Given $\epsilon > 0$, it is enough to show that the singularities of $Pu$ lie within an $\epsilon$ neighbourhood of those of $u$. To see this we decompose both $K$ and $u$ into pieces supported near the singularities and pieces which are smooth. So suppose

$$K = K_1 + K_2$$

with the support of $K_1$ contained in

$$\{|x - y| < \epsilon/2\}$$

and $K_2$ smooth and

$$u = u_1 + u_2$$

with $\text{supp}(u_1)$ within an $\epsilon/2$ neighbourhood of $\text{singsupp}(u_1)$ and $u_2$ smooth. Then

$$Ku = K_1u_1 + K_2u_1 + Ku_2.$$

Now from the relation of supports, we have

$$\operatorname{supp}(K_1 u_1) \subset \{x : |x - y| < \epsilon, y \in \operatorname{singsupp} u\}$$

and since the singular support is smaller than the support, the first term is OK. We also have $Ku_2$ is smooth as $u_2$ is.

This leaves $K_2 u_1$. This turns out to be smooth also as $K_2$ is and the result will then follow. This is an important result so we will make it a theorem. $\qquad\square$

**Theorem 4.6.** *If $K \in C^\infty(\mathbb{R}^n \times \mathbb{R}^n)$ and is properly supported then $K$ induces a map from $\mathcal{D}'(\mathbb{R}^n)$ to $C^\infty(\mathbb{R}^n)$.*

*Proof.* Consider the map

$$Ku : x \mapsto < K(x, y), u(y) > .$$

This is well-defined for each $x$, as $K(x, y)$ is a smooth function in $y$. We want to show the result is smooth in $x$.

But by the definition of a distribution, we have for some $k$,

$$\|Ku(x) - Ku(x')\| \leq \sum_{|\alpha| \leq k} \sup_y |D_y^\alpha (K(x, y) - K(x', y))|.$$

The properness of $K$ will ensure that the supremum is over a compact set and thus the function will certainly be continuous.

Note the same is true for the maps

$$D_x^\alpha : x \mapsto < D_x^\alpha K(x, y), u(y) >$$

so the formal derivatives are continuous too. So, as in the proof that pseudo-differential operators are continuous, $Ku$ is smooth as required. $\qquad\square$

Since we have shown that smooth kernels kill all singularities, we will ignore them most of the time. So we shall regard two operators as equivalent if and only their Schwartz kernels differ by a smooth term and we write this as

$$P \equiv Q.$$

We shall call an inverse up to a smoothing term a *parametrix*. Our theorem shows that such an operator will allow us to describe the singularities of $u$ from those of $Pu$. Note that if $K$ is smoothing then $PK, KP$ are smoothing also for any pseudo-differential operator so the smoothing operators form an ideal in the class of pseudo-differential operators.

The next thing we want to do is prove that the composition of two pseudo-differential operators is a pseudo-differential operator. In order to do this, we need a better understanding of the relationship between the symbol $a$ and the associated operator $P$.

**Theorem 4.7.** *If $P$ is a pseudo-differential operator of order $m$, then up to smooth terms the Schwartz kernel of $P$ can be written in the form*

$$\int e^{i\langle x-y,\xi\rangle} c(x,\xi) d\xi$$

*or in the form*

$$\int e^{i\langle x-y,\xi\rangle} d(y,\xi) d\xi$$

*with $c, d \in S^m(\mathbb{R}^n; \mathbb{R}^n)$.*

Note these are called the left and right quantizations.

*Proof.* So we have initially that

$$K(x,y) = \int e^{i\langle x-y,\xi\rangle} a(x,y,\xi) d\xi.$$

We want to eliminate the dependence on one of the base variables. Recall from Taylor's theorem that

$$a(x,y,\xi) = \sum_{|\alpha| \leq N-1} \frac{1}{\alpha!} \partial_y^\alpha a(x,x,\xi)(y-x)^\alpha + \sum_{|\alpha|=N} (y-x)^\alpha b_\alpha(x,y,\xi)$$

and $b_\alpha = \int_0^1 (1-t)^{N-1} \partial_y^\alpha a(x,(1-t)x+ty,\xi) dt$. It is immediate (and important) that $\partial_y^\alpha a(x,x,\xi) \in S^m$ and $b_\alpha(x,y,\xi) \in S^m$.

Let's consider the terms $\partial_y^\alpha a(x,x,\xi)(y-x)^\alpha$,

$$\int e^{i\langle x-y,\xi\rangle} (x-y)^\alpha \partial_y^\alpha a(x,x,\xi) d\xi = \int D_\xi^\alpha e^{i\langle x-y,\xi\rangle} \partial_y^\alpha a(x,x,\xi) d\xi$$

and so integrating by parts we obtain

$$\int e^{i\langle x-y,\xi\rangle} (x-y)^\alpha \partial_y^\alpha a(x,x,\xi) d\xi = \int e^{i\langle x-y,\xi\rangle} (-1)^{|\alpha|} D_\xi^\alpha \partial_y^\alpha a(x,x,\xi) d\xi.$$

(check this actually works!)

We have applying the same argument to $b_\alpha$ that the result is true up to a term of order $m - N$. To proceed further, we need to sum the series

$$\sum_\alpha \frac{1}{\alpha!} D_\xi^\alpha \partial_y^\alpha a(x, x, \xi).$$

This sum will not converge but we can make sense of it in an asymptotic way. We will use the following result:

**Proposition 4.8.** *If $a_j \in S^{m-j}$ then there exists $a \in S^m$ such that*

$$a - \sum_{j=0}^{N-1} a_{m-j} \in S^{m-N}, \ \forall N$$

The proof of this is on the example sheet 1. Note an immediate corollary of this, which will come in handy later on,

**Corollary 4.9.** *If $P_j \in \Psi DO^{m-j}(\mathbb{R}^n)$ for $j = 0, \ldots, \infty$ then there exists $P \in \Psi DO_{cl}^{m-j}(\mathbb{R}^n)$ such that*

$$P - \sum_{j<N} P_j \in \Psi DO^{m-N}(\mathbb{R}^n)$$

*for all $N$.*

So with this proposition, we let $c$ be an asymptotic sum of

$$\sum_\alpha \frac{1}{\alpha!} D_\xi^\alpha \partial_y^\alpha a(x, x, \xi)$$

and is then clear that

$$\int e^{i\langle x-y, \xi \rangle} c(x, \xi) d\xi - \int e^{i\langle x-y, \xi \rangle} a(x, y, \xi) d\xi \in \Psi DO^{m-N}, \ \forall N$$

which establishes the first half of the result.

The second half is identical, except that the Taylor expansion is in $y$ not $x$.                                                                $\square$

Note we have actually done better than the statement of our theorem in that the proof gives us asymptotic formulaes for $c, d$

$$c(x, \xi) \sim \sum_\alpha \frac{1}{\alpha!} D_\xi^\alpha \partial_y^\alpha a(x, x, \xi), \ d(y, \xi) \sim \sum_\alpha \frac{(-1)^{|\alpha|}}{\alpha!} D_\xi^\alpha \partial_x^\alpha a(y, y, \xi).$$

Note that unlike our original symbol which depended on $x, y$ and was definitely not unique, these symbols are essentially unique.

If

$$K(x, y) = \left(\frac{1}{2\pi}\right)^n \int e^{i\langle x-y, \xi\rangle} a(x, \xi) d\xi$$

then this is really just a Fourier transform across $x = y$. Putting $w = x - y, x' = x$, we have

$$L(x', w) = K(x', x' - w) = \left(\frac{1}{2\pi}\right)^n \int e^{i<w,\xi>} a(x', \xi) d\xi$$

and so

$$a(x', \xi) = \int e^{-iw.\xi} L(x', w) dw.$$

i.e. $a$ is determined by $K$. A smooth error on $K$, for $K$ properly supported, would translate into an error of order $-\infty$ on $a$, as the Fourier transform of a compactly supported smooth function is Schwartz.

**Definition 4.10.** *If $P$ has Schwartz kernel*

$$\int e^{i\langle x-y, \xi\rangle} a(x, \xi) d\xi$$

*with $a \in S^m$, up to smooth terms, then the total left symbol of $P$, $\sigma_L(P)$, is $a(x, \xi)$. We will regard $\sigma_L(P)$ as an element of $S^m/S^{-\infty}$.*

So the total left symbol determines $P$ up to order $-\infty$.

The same arguments work on the other side and we make a similar definition.

**Definition 4.11.** *If $P$ has Schwartz kernel*

$$\int e^{i\langle x-y, \xi\rangle} b(y, \xi) d\xi$$

*with $b \in S^m$ ,up to smooth terms, then the total right symbol of $P$, $\sigma_R(P)$, is $b(y, \xi)$. We will regard $\sigma_R(P)$ as an element of $S^m/S^{-\infty}$.*

Why use left and right symbols? If we have operators $P_j \in \Psi DO^{m_j}(\mathbb{R}^n)$ then giving $P_1$ the left form and $P_2$ the right one, they can be written (up to smooth terms)

$$P_1 u = \left(\frac{1}{2\pi}\right)^n \int e^{ix.\xi} p(x, \xi) \hat{u}(\xi) d\xi$$

and

$$P_2 v = \left(\frac{1}{2\pi}\right)^n \int e^{ix.\xi} \int q(y, \xi) v(y) dy d\xi.$$

So it is then immediate from the Fourier inversion theorem that

$$P_1 P_2 v = \left(\frac{1}{2\pi}\right)^n \int e^{i\langle x-y,\xi\rangle} p(x,\xi) q(y,\xi) u(y) dy d\xi.$$

So $P_1 P_2$ is a pseudo-differential operator of order $m_1 + m_2$ with total symbol $p(x,\xi) q(y,\xi)$. Of course, we want a more symmetric expression in terms of outputs and inputs - we want to compute the left symbol of the composite in terms of the left symbols of the inputs.

**Theorem 4.12.** *If $P_j \in \Psi DO^{m_j}$ are properly supported and have total symbols $p_j(x,\xi)$ then $P_1 P_2 \in \Psi DO^{m_1+m_2}$ and has total symbol*

$$a(x,\xi) \sim \sum_\alpha \frac{i^{|\alpha|}}{\alpha!} D_\xi^\alpha p_1(x,\xi) D_x^\alpha p_2(x,\xi).$$

The only thing left to prove is the correctness of the asymptotic expansion. This is just an exercise in binomial expansions - first convert the left symbol of $p_2$ to a right one and then convert the product to a left one.

## 5. Classicality

So far we have been working with general pseudo-differential operators - this class is in fact larger than we need and we can work with a smaller class which has some nice properties and will include all the operators we actually need. Recall that a homogeneous function smoothed off at the origin is a symbol so we can require a symbol to have an asymptotic expansion in homogeneous terms.

**Definition 5.1.** *A symbol $a(x,\theta)$ of order $m$ is classical if there exists a sequence of functions $a_{m-j}(x,\theta)$, homogeneous of order $m - j$ in $\theta$ such that*

$$a - \sum_{j=0}^{N-1} (1-\phi)(\theta) a_{m-j}(x,\theta) \in S^{m-N}$$

*with $\phi$ a compactly supported function equal to 1 near 0.*

**Definition 5.2.** *A pseudo-differential operator is classical if its total left (or right) symbol is classical. We denote this class $\Psi DO_{cl}^m(\mathbb{R}^n)$*

Note that the asymptotic expansion is unique and so knowing the operator and knowing the expansion of its symbol are really the same thing up to smooth terms which we always ignore!

It is obvious that composition preserves classicality and that differential operators are classical pseudo-differential operators. The most important thing about classical operators is that they have a unique principal symbol.

**Definition 5.3.** *If $P \in \Psi DO_{cl}^m(\mathbb{R}^n)$ then if $P$ has Schwartz kernel*

$$\int e^{i\langle x-y,\xi\rangle} p(x,\xi) d\xi$$

*with asymptotic expansion*

$$\sum p_{m-j}(x,\xi)$$

*then the principal symbol of $P$ is $p_m(x,\xi)$ and is denoted $\sigma_m(P)(x,\xi)$.*

We denote the space of smooth homogeneous functions on $\mathbb{R}_x^n \times (\mathbb{R}_\theta^n - \{0\})$ which are homogenous of degree $m$ in $\theta$ by $H^m(\mathbb{R}^n; \mathbb{R}^n - 0)$. We express the fact that $\sigma_m(P)(x,\xi)$ determines $P$ up to one order lower by using a short exact sequence. All this means is that the image of one map is exactly equal to the kernel (null-space) of the next.

**Theorem 5.4.** *The following sequence is exact*

$$0 \longrightarrow \Psi DO_{cl}^{m-1}(\mathbb{R}^n) \longrightarrow \Psi DO_{cl}^m(\mathbb{R}^n) \xrightarrow{\sigma_m} H^m(\mathbb{R}^n; \mathbb{R}^n - 0) \longrightarrow 0.$$

The reason the principal symbol is particularly useful is that the formula for products is nice.

**Theorem 5.5.** *If $P \in \Psi DO_{cl}^m(\mathbb{R}^n)$, $Q \in \Psi DO_{cl}^{m'}(\mathbb{R}^n)$ are properly supported and have classical expansions*

$$\sum p_{m-j}(x,\xi), \sum q_{m'-k}(x,\xi)$$

*then $PQ$ is also classical with asymptotic expansion $\sum r_{m+m'-l}(x,\xi)$ where*

$$r_{m+m'-l} = \sum_{j+k+|\alpha|=l} \frac{1}{\alpha!} D_\xi^\alpha p_{m-j}.\partial_x^\alpha q_{m'-k}$$

*and thus*

$$\sigma_{m+m'}(PQ) = \sigma_m(P)\sigma_{m'}(Q).$$

*Proof.* We know $PQ$ is a pseudo-differential operator with asymptotic expansion

$$\sum_\alpha \frac{1}{\alpha!} D_\xi^\alpha p(x,\xi) \partial_x^\alpha q(x,\xi).$$

Now $D_\xi^\alpha p(x,\xi)$ has asymptotic expansion

$$\sum D_\xi^\alpha p_{m-j}$$

and $D_\xi^\alpha p_{m-j}$ is homogeneous of degree $m-j-|\alpha|$ and $\partial_x^\alpha q(x,\xi)$ has asymptotic expansion

$$\sum \partial_x^\alpha q_{m'-k}$$

with $\partial_x^\alpha q_{m'-k}$ homogeneous of degree $m'-k$. The result follows by collecting the terms at each level of homogeneity.

The statement about principal symbols now follows from considering the top level of homogeneity. □

Now we have done all this, we can start constructing parametrices - ie inverses up to terms of order $-\infty$.

**Definition 5.6.** *A classical pseudo-differential operator is said to be elliptic of order $m$ if $P \in \Psi DO_{cl}^m$ and $\sigma_m(P)$ is never zero.*

**Example 5.7.** The Laplacian

$$\Delta = \sum D_{x_j}^2$$

plus any first order perturbation has principal symbol, $\sum\limits_{j=1}^n \xi_j^2$, and is therefore elliptic. Or more generally, if $g_{ij}(x)$ is a Riemannian metric on $\mathbb{R}^n$ and $g^{ij}(x)$ is the inverse matrix then the variable coefficient Laplacian

$$\Delta = \sum g^{ij} D_{x_i} D_{x_j} + E,$$

with $E$ first order, has principal symbol $\sum g^{ij}\xi_i\xi_j$, and is therefore elliptic. On $\mathbb{R}^2$, we also have that the Neumann operator

$$D_x + i D_y$$

has principal symbol $\xi + i\eta$ and is therefore elliptic.

**Theorem 5.8.** *If $P$ is a properly supported, elliptic, classical pseudo-differential operator of order $m$ then there exists*

$$Q \in \Psi DO_{cl}^{-m}$$

*such that*

$$PQ - \mathrm{Id}, QP - \mathrm{Id} \in C^\infty.$$

*Proof.* First, note that the principal and total symbol of Id are both 1. So we know that for the result to be true

$$\sigma_m(P)\sigma_{-m}(Q) = 1.$$

So we put $\sigma_{-m}(Q) = \sigma_m(P)^{-1}$. The assumption of ellipticity ensures that $\sigma_{-m}(Q)$ is smooth and thus we can pick such a $Q$.

Now our formula for products of symbols tells us that the $l^{th}$ term in the expansion for $PQ$ will be

$$\sum_{|\alpha|+j+k=l} \frac{1}{\alpha!} D_\xi^\alpha p_{m-j}(x,\xi) \partial_x^\alpha q_{m-k}(x,\xi) \qquad (3)$$

and we need this to be equal to zero for each $l \geq 1$.

So suppose that $q_{-m}, \ldots, q_{-m-N}$ have been chosen so that (3) is satisified for $l \leq N$ then we can define

$$q_{-m-N-1} = -p_m^{-1} \sum_{|\alpha|+j+k=N+1, k<N+1} \frac{1}{\alpha!} D_\xi^\alpha p_{m-j}(x,\xi) \partial_x^\alpha q_{-m-k}$$

and (3) will be satisfied for $l = N+1$ also. Thus by induction, we can construct $q_{-m-k}$ so that (3) is true for all $k$.

Let $q$ have classical expansion $q_{m-k}$ and then we have

$$PQ - \mathrm{Id} \in \Psi DO^{-\infty} = C^\infty.$$

A similar argument constructs $Q'$ such that

$$Q'P - \mathrm{Id} \in \Psi DO^{-\infty} = C^\infty.$$

But regarding $\Psi DO_{cl}^m / \Psi DO^{-\infty}$ as a semi-group, left inverses and right inverses must be equal and the result follows. $\qquad\square$

We can now prove elliptic regularity.

**Theorem 5.9.** *(Weyl's Lemma) If $P \in \Psi DO^m$ is proper, classical and $P$ is elliptic of order $m$ and $u \in \mathcal{D}'(\mathbb{R}^n)$ and*

$$Pu \in C^\infty$$

*then*

$$u \in C^\infty.$$

*More generally, for any $u \in \mathcal{D}'(\mathbb{R}^n)$,*

$$\mathrm{singsupp}(Pu) = \mathrm{singsupp}(u).$$

*Proof.* Let $Q$ be a parametrix for $P$ then $QP - \mathrm{Id} = R$ and $Ru \in C^\infty$.

So by pseudo-locality, theorem 4.5,

$$\begin{aligned}
\mathrm{singsupp}(u) &= \mathrm{singsupp}((\mathrm{Id}+R)u) \\
&= \mathrm{singsupp}(QPu) \\
&\subset \mathrm{singsupp}(Pu) \\
&\subset \mathrm{singsupp}(u).
\end{aligned}$$

□

Note this is actually a deep result and is quite hard to prove by other means! It also relies heavily on the fact that $P$ is elliptic and is definitely not true for a general pseudo-differential operator.

## 6. Continuity

So, we have seen that elliptic operators have parametrices which provide inverses up to smooth terms. i.e. they provide an inverse at the singularity level. Sobolev Spaces can be used to measure the badness of a singularity and we will examine how a pseudo-differential operator maps between them.

**Definition 6.1.** *We shall say $u \in H^s(\mathbb{R}^n)$ if $u \in \mathcal{S}'(\mathbb{R}^n)$ and*

$$\hat{u}(\xi)(1 + |\xi|^2)^{s/2} \in L^2(\mathbb{R}^n).$$

Note that the lower $s$ is, the more singular $u$ becomes.

It is trivial that

$$D_{x_j} : H^s \to H^{s-1}.$$

There are two ways to fail to be in $L^2$ - a distribution can be too big at $\infty$ or too singular at some point. The same is true of $H^s$ so, since what we really want to measure is how singular a distribution is, we will often work with two related spaces.

**Definition 6.2.**

$$H_c^s(\mathbb{R}^n) = \{u : u \in H^s(\mathbb{R}^n) \ and \ supp(u) \ compact.\}.$$

$$H_{loc}^s(\mathbb{R}^n) = \{u : \phi u \in H^s(\mathbb{R}^n), \forall \phi \in C_0^\infty\}$$

*These are the compact and local spaces.*

The theorem we want to prove is

**Theorem 6.3.** *If $P \in \Psi DO_{cl}^m(\mathbb{R}^n)$ and is properly supported then*

$$P : H_c^s \to H_c^{s-m}.$$

Most of the work involved is in proving the simplest case - $m = 0, s = 0$. ie the $L_c^2$ continuity of zeroth order pseudo-differential operators.

To do this we will need to understand adjoints of pseudo-differential operators - we already understand transposes and these are similar.

Recall, if $\phi, \psi \in C_0^\infty$ then the transpose of $P$, is defined by

$$\int P\phi(x)\psi(x)dx = \int \phi(x)P^t\psi(x)dx$$

whereas the adjoint is defined by

$$\int P\phi(x)\overline{\psi(x)}dx = \int \phi(x)\overline{P^*\psi(x)}dx.$$

So these operators are pretty similar and indeed are related by

$$P^*\psi = \overline{P^t\bar{\psi}}.$$

Thus, if $P \in \Psi DO^m$ has total left symbol $a(x, \xi)$ we deduce that $P^*$ is pseudo-differential and has a total (right) symbol

$$\bar{a}(y, \xi)$$

and so has total left symbol $b(x, \xi)$ with asymptotic expansion

$$\sum_\alpha \frac{1}{\alpha!} \partial_x^\alpha D_\xi^\alpha \bar{a}(x, \xi).$$

So $P^*$ will be classical if $P$ is, and

$$\sigma_m(P^*) = \overline{\sigma_m(P)}.$$

Thus we obtain

**Theorem 6.4.** *If $P \in \Psi DO_{cl}^m$ is self-adjoint i.e. $P = P^*$ then $\sigma_m(P)$ is real-valued.*

*If $f(x, \xi)$ is homogeneous of order $m$ in $\xi$ then there exists a pseudo-differential operator which is self-adjoint with principal symbol $f$ if and only if $f$ is real-valued.*

*Proof.* The only part of this which is not wholly obvious is that if $f$ is real-valued then there exists a self-adjoint operator with principal symbol $f$. To see this, let $\sigma_m(Q) = f$ some $Q \in \Psi DO_{cl}^m(\mathbb{R}^n)$ and put

$$P = \frac{Q + Q^*}{2}.$$

$\square$

We prove continuity using an approximate square root which will reduce to proving the result in the case where $P$ has smooth kernel.

**Theorem 6.5.** *If $P \in \Psi DO_{cl}^m$ is self-adjoint and $\sigma_m(P)$ is positive then there exists $Q \in \Psi DO_{cl}^{m/2}$ which is self-adjoint such that*

$$P - Q^2 \in \Psi DO^{-\infty}.$$

*Proof.* For a first guess we can take $Q_0 \in \Psi DO_{cl}^{m/2}$ self-adjoint with principal symbol $\sigma_m(P)^{1/2}$ which is a valid symbol by our hypothesis.

Then we have

$$P - Q_0^2 \in \Psi DO^{m-1}.$$

We use a recursive argument so now suppose we have constructed $Q_0, Q_1, \ldots, Q_{N-1}$ classical and self-adjoint such that

$$P - \left(\sum_{j=0}^{N-1} Q_j\right)^2 = R_N \in \Psi DO_{cl}^{m-N}$$

and we want to construct $Q_N \in \Psi DO_{cl}^{m/2-N}$.

Now

$$P - \left(\sum_{j=0}^{N} Q_j\right)^2 = R_N - Q_N \left(\sum_{j=0}^{N-1} Q_j\right) - \left(\sum_{j=0}^{N-1} Q_j\right) Q_N - Q_N^2.$$

The principal symbol of this will be

$$\sigma_{m-N}(R_N) - 2\sigma_{m/2-N}(Q_N)\sigma_{m/2}(Q_0)$$

and to continue the recursion this has to equal zero. So noting that $R_N$ is self-adjoint we can pick $Q_N$ to be self-adjoint with principal symbol equal to

$$\frac{1}{2}\sigma_{m/2}(Q_0)^{-1}\sigma_{m-N}(R_N)$$

which is of course real.

So once we have $Q_j$ constructed for all $j$ then we can pick $Q'$ such that

$$Q' - \sum_{j<N} Q_j \in \Psi DO^{m-N}, \ \forall N.$$

To do this just asymptotically sum at the symbol level. $Q'$ will necessarily be classical as all the $Q_j$ are. However $Q'$ is not obviously self adjoint so we put $Q = \frac{Q'+Q'^*}{2}$ and we have constructed the requisite $Q$. Note that $Q - Q'$ is of order $-\infty$ as each $Q_j$ is self-adjoint. So $Q$ is also an asymptotic sum for the $Q_j$.

By construction $Q$ will now be the approximate square root. $\square$

Proving $L_c^2$ continuity of a properly-supported, zeroth-order pseudo-differential operator, $A$, means that for each compact set $K$ we have to show for $\phi \in L^2, \text{supp}(\phi) \subset K$ that there exists $K'$ compact and $C_K$ such that

$$\|A\phi\| \leq C_K \|\phi\|$$

and supp$(A\phi) \subset K'$ with the norms taken into the $L^2$ sense. In fact as compactly supported smooth functions are dense in $L_c^2$ we need only prove this for smooth functions. Note that if we only work with smooth functions then we establish the existence of continuous extension, $\tilde{A}$, to $L_c^2$. But this will agree with $A$ on $L_c^2$ because $A$ is weakly continuous on distributions and agrees with $\tilde{A}$ on $C_0^\infty$.

The existence of a $K'$ as above is immediate from the relation of supports so we need to establish the bound $C_K$. Since only the behaviour of the kernel of $A$ on the set $\pi_y^{-1}(K)$ will affect its mapping properties on elements of $L^2(K)$, we can decompose the kernel of $A$ into one piece which is compactly supported and another which maps all elements of $L^2(K)$ to zero. Thus it is enough to consider the case supp$(A)$ compact.

It is easier to work with squares and we have

$$\|A\phi\|^2 = \langle A\phi, A\phi \rangle$$
$$= \langle \phi, A^*A\phi \rangle.$$

Now $A^*A$ is self-adjoint and its principal symbol is $|\sigma_0(A)|^2$ which is necessarily non-negative but may have some zeroes. The fact that the kernel of $A$ is compactly supported means that $|\sigma_0(A)|^2$ will be compactly supported on the base and since it is homogeneous of degree zero, we deduce from continuity that it must be bounded above. So suppose $|\sigma_0(A)|^2 < C$ everywhere.

Then the operator $C - A^*A$ is self-adjoint and has strictly positive principal symbol, thus by our theorem 6.5 we have

$$C - A^*A = B^2 + G$$

with $B$ self-adjoint and $G$ smooth (and self-adjoint).

So we have,

$$\|A\phi\|^2 = C\|\phi\|^2 - \langle \phi, B^2\phi \rangle - \langle \phi, G\phi \rangle.$$

The second term equals $-\|B\phi\|^2$ and thus is necessarily negative so it remains simply to show that

$$|\langle \phi, G\phi \rangle| \leq C'\|\phi\|^2.$$

Using Cauchy-Schwartz, it is enough to show the $L_c^2$ continuity of $G$.

As before, is enough to consider $G$ compactly supported and since $G$ is smooth, it must also be bounded (as a function) by say $D$.

So we have,

$$\|G\phi\|^2 = \int \left| \int G(x,y)\phi(y)dy \right|^2 dx$$

and since $|G(x,y)| \leq D$ and our integrals are supported in the finite domain $K \times K'$, we deduce

$$\|G\phi\|^2 \leq D^2 V \|\phi\|^2$$

where $V$ is volume of $K \times K'$.

Thus we have proven,

**Theorem 6.6.** *If $A \in \Psi DO_{cl}^0(\mathbb{R}^n)$ is classical and properly supported then $A$ is continuous as a map from $L_c^2$ to $L_c^2$.*

The fact we are considering only compactly supported functions is essential here - otherwise consider for example multiplication by $e^x$ - this is clearly not continuous on $L^2$ but is a zeroth order differential operator.

We now wish to prove in general that properly supported elements of $\Psi DO_{cl}^m(\mathbb{R}^n)$ map

$$H_c^s \to H_c^{s-m},$$

we do this by reducing to the zeroth order case.

So let $\Lambda_s$ be a properly supported element of $\Psi DO^s$ with total symbol $\langle \xi \rangle^s$. This is classical as

$$\langle \xi \rangle^s = (1 + |\xi|^2)^{s/2} = |\xi|^s(1 + |\xi|^{-2})^{s/2},$$

and we can apply a binomial expansion.

We will show that $\Lambda_s$ defines a continuous map from $H_c^t$ to $H_c^{t-s}$. Then if $P \in \Psi DO_{cl}^m$ is properly supported, we have that

$$P = \Lambda_s \left( \Lambda_{-s} P \Lambda_{s-m} \right) \Lambda_{m-s} + R$$

with $R$ of order $-\infty$. Now, $(\Lambda_{-s} P \Lambda_{s-m}) \in \Psi DO_{cl}^0$ and is properly supported, so it defines a continuous map of $L_c^2$. So if we know that $R$ is continuous, it immediately follows that $P$ is continuous as a map from

$$H_c^s \to H_c^{s-m}.$$

(Recall that $H^0 = L^2$.)

We can write

$$\Lambda_s u = M_s u + R' u$$

where

$$M_s u = \left( \frac{1}{2\pi} \right)^n \int e^{ix.\xi} \langle \xi \rangle^s \hat{u}(\xi) d\xi$$

where $R'$ is of order $-\infty$. Unfortunately, $M_s$ will not be properly supported but if $\text{supp}(u) \subset K$, $K$ compact, we have that $\text{supp}(\Lambda_s u)$ is contained in

$$K' = \text{supp}(\Lambda_s) \circ K.$$

So if we take $\psi = 1$, $\psi \in C_0^\infty(\mathbb{R}^n)$ on $K'$, we have that

$$\Lambda_s u = \psi M_s \psi u + \psi R \psi u$$

and $\psi R \psi$ will have kernel $\psi(x) R(x,y) \psi(y)$.

Thus to show the continuity of $\Lambda_s$, it is enough to show three things

1. $M_s : H^t \to H^{t-s}$
2. Multiplication by $\psi$ is continuous map from $H^s$ to $H_c^s$.
3. If $R \in \Psi DO^{-\infty}$ and is compactly supported then $R$ defines a continuous map from $H_c^s \to H_c^N$ for any $N$.

The most innocuous looking of these three, 2., is in fact the hardest! Let's start with 1. Applying $M_s$ is equivalent to multiplying on the Fourier transform side by

$$\langle \xi \rangle^s,$$

which by definition is an isometry from $H^t \to H^{t-s}$.

To see 3., apply the Fourier transform to obtain

$$\widehat{Ru}(\xi) = \int \tilde{r}(\xi, y) u(y) dy.$$

Here $\tilde{r}$ is the partial Fourier transform. We thus deduce,

$$\widehat{Ru}(\xi) = (2\pi)^n \int \hat{r}(\xi, \eta) \hat{u}(-\eta) d\eta.$$

Of course, $\hat{r}(\xi, \eta)$ is Schwartz and so is less than or equal to

$$C \langle \xi \rangle^{-M_1} \langle \eta \rangle^{-M_2}$$

for any $M_1, M_2$. So we conclude,

$$|\widehat{Ru}(\xi)| \langle \xi \rangle^{2N} \leq \left( \int \langle \xi \rangle^N |\hat{r}(\xi, \eta)|^2 \langle \eta \rangle^{-2s} d\eta \right)^{\frac{1}{2}} \|u\|_s.$$

The integral will be rapidly decreasing in $\xi$ and so we have that $|\hat{Ru}(\xi)| \langle \xi \rangle^N \in L^2$ with norm less than or equal to $C\|u\|_s$ and the continuity follows.

This leaves us with the problem of proving 2. We first prove this for $s \geq 0$ and then use *duality* to deduce the general case.

Let $M_s$ be the operator defined by

$$\widehat{M_s u}(\xi) = \langle \xi \rangle^s \hat{u}(\xi)$$

this defines an isometry from

$$H^t \to H^{t-s}$$

for each $t$. So $M_\psi$ is a continuous map on $H^s$ if and only if

$$T_s = M_s M_\psi M_{-s}$$

is a continuous map on $H^0 = L^2$. Since $C_c^\infty$, is dense in $L^2$ we need only check this on $C_c^\infty$. So let $u \in C_c^\infty$ and then $\widehat{T_s u}$ is given by

$$\langle \xi \rangle^s \int \hat{\psi}(\xi - \xi') < \xi' >^{-s} \hat{u}(\xi') d\xi'.$$

Using the fact that

$$(\xi_j - \xi_j') \hat{\psi}(\xi - \xi') = \widehat{D_j \psi}(\xi - \xi'),$$

we can rewrite this for even $N$ as

$$\langle \xi \rangle^s \int \langle \xi - \xi' \rangle^{-N} \widehat{\psi_N}(\xi - \xi') < \xi' >^{-s} \hat{u}(\xi') d\xi',$$

for some $\psi_N \in C_0^\infty(\mathbb{R}^n)$. Noting that

$$\langle \xi \rangle \leq C < \xi' > \langle \xi - \xi' \rangle,$$

for $s \geq 0$ we have

$$|\widehat{T_s u}| \leq C \int |\widehat{\psi_N}(\xi - \xi')||\hat{u}(\xi')| d\xi'.$$

Recall that the Fourier transform is an isometry of $L^2$ so it is enough to estimate the $L^2$ norm of $|\widehat{T_s u}|$. Now

$$\int |\widehat{T_s u}(\xi)|^2 d\xi \leq \int \left( \int |\widehat{\psi_N}(\xi - \xi')||\hat{u}(\xi')| d\xi' \right)^2 d\xi$$

$$= \int \left( \int |\widehat{\psi_N}(\xi - \xi')|^{\frac{1}{2}} |\widehat{\psi_N}(\xi - \xi')|^{\frac{1}{2}} |\hat{u}(\xi')| d\xi' \right)^2 d\xi$$

so applying Cauchy-Schwartz,

$$\leq \int \left( \int |\widehat{\psi_N}(\xi - \xi')||\hat{u}(\xi')|^2 d\xi' \int |\widehat{\psi_N}(\xi - \xi')| d\xi' \right) d\xi$$

$$= \int \int |\widehat{\psi_N}(\xi - \xi')||\hat{u}(\xi')|^2 d\xi' d\xi \|\widehat{\psi_N}\|_{L_1}$$

putting $\eta = \xi - \xi'$, we have

$$= \int |\hat{u}(\xi')|^2 d\xi' \|\widehat{\psi_N}\|_{L_1}^2.$$

Thus we have shown,

$$\|T_s u\|_2 \leq \|\widehat{\psi_N}\|_1 \|u\|_2$$

and so $M_\psi$ is continuous for $s \geq 0$. The technique we will use to prove the result for general $s$ is sufficiently important that we will give it its own section.

## 7. Duality For Sobolev Spaces

There is a natural pairing between $H^s$ and $H^{-s}$, we define

$$\langle .,. \rangle : H^s \times H^{-s} \to \mathbb{C} \tag{4}$$

$$\langle u, v \rangle = \int \hat{u}(\xi)\overline{\hat{v}(\xi)}d\xi. \tag{5}$$

We can think of this as the $L^2$ pairing of $\langle \xi \rangle^{-s}\hat{u}$ and $\langle \xi \rangle^s\hat{v}$ which shows that the integral makes sense. This pairing is linear in the first factor and quasi-linear in the second. It is also continuous - using Cauchy-Schwartz for the $L^2$ pairing, we have $|\langle u, v \rangle| \leq \|u\|_s\|v\|_{-s}$.

This means that the pairing gives a natural map from $H^s$ to the dual space of $H^{-s}$. Simply, define

$$T : H^s \to (H^{-s})^*, \tag{6}$$

$$Tv(u) = <v, u> \tag{7}$$

The continuity of the pairing yields that $Tv$ is in $(H^{-s})^*$ and that the norm of $Tv$ is less than or equal to $\|v\|_s$.

Now if $v \in H^s$ then we can define $u$ in $H^{-s}$ via

$$\hat{u}(\xi) = \langle \xi \rangle^{2s}\hat{v}(\xi)$$

and then

$$|Tv(u)| = Tv(u) = \int |\hat{v}(\xi)|^2\langle \xi \rangle^{2s}d\xi = \|v\|_s^2$$

and so

$$\|Tv\| = \|v\|_s.$$

So $T$ is a quasi-linear isometry. It is in fact surjective. Recall that as $H^{-s}$ is a Hilbert space, it is naturally isomorphic to $(H^{-s})^*$. So given an element, $f$, of $(H^{-s})^*$, we have from the Riesz representation theorem an element $v \in H^{-s}$ such that

$$f(u) = \langle u, v \rangle_s, \forall u \in H^{-s}.$$

But

$$\langle u, v \rangle_{-s} = \int \hat{u}(\xi)\overline{\hat{v}(\xi)}\langle \xi \rangle^{-2s}d\xi.$$

So taking $v_1$ in $H^s$ such that $\hat{v_1}(\xi) = \langle \xi \rangle^{-2s}\hat{v}(\xi)$, we have that

$$Tv_1 = f.$$

So $T$ is onto.

Recall that we wanted to show that multiplication by a compactly supported smooth function is continuous on Sobolev spaces and we had already shown this for $s \geq 0$.

The adjoint of multiplication by $\psi$ on $H^s$ is just multiplication by $\bar{\psi}$ on $H^{-s}$.

So for $s < 0$, $u \in H^{-s}, v \in H^s$ we deduce,

$$
\begin{aligned}
|\langle M_\psi v, u \rangle| &= |\langle v, M_{\overline{\psi}} u \rangle| \\
&\leq \|v\|_s \|M_{\overline{\psi}} u\|_{-s} \\
&\leq C\|v\|_s \|u\|_{-s}.
\end{aligned}
$$

This says precisely that as an operator the norm of $M_\psi v$ is less than or equal to $C\|v\|_s$. As the operator norm equals the original norm, we deduce that $M_\psi$ is continuous. Note that the idea of this proof is very important and is a standard technique.

## 8. Coordinate Invariance

We know that under changes of coordinates differential operators remain differential operators and we can compute the new operator by repeated use of the chain rule. We want to do the same thing for pseudo-differential operators - ie prove they remain pseudo-differential and compute the symbol of the resulting operator. It turns out that the lower order terms of the symbol are very complicated but that the principal symbol transforms very nicely.

What do we really mean by a change of coordinates, well suppose we have on operator on an open subset $U$ of $\mathbb{R}^n$ and the operator $P$ is a map from compactly supported smooth functions on $U$ to smooth functions on open $V$ containing $U$. Then if we have another pair of sets $U_1, V_1$ and a map $\chi$ from $V_1$ to $V$ such that $\chi(U_1) = U$ and $\chi$ is a smooth bijection and has a smooth inverse then we can define an operator from compactly supported smooth functions on $U_1$ to smooth functions on $V_1$ by

$$
P_1 f = \left[ P(f \circ \chi^{-1}) \right] \circ \chi.
$$

That is we take the function $f$ on $U_1$ and make it a function on $U$ by pulling back by $\chi^{-1}$ i.e. taking the map

$$
\left( \chi^{-1} \right)^* f : U \to \mathbb{C}
$$

defined by

$$\left(\chi^{-1}\right)^* f(x) = f(\chi^{-1}x)$$

and then we apply $P$ to get an function on and then by we pull-back the resulting functon by $\chi$ to get a function on $V_1$. Thus we have a map, $P_1$, from $C_0^\infty(U_1) \rightarrow C^\infty(V_1)$.

So now suppose $P$ is a pseudo-differential operator, we can write

$$P_1 u(x) = \int e^{i<\chi(x)-y,\xi>} a(\chi(x), y, \xi) u(\chi^{-1}(y)) dy d\xi.$$

Doing a change of coordinates, $y = \chi(z)$ this becomes,

$$P_1 u(x) = \int e^{i<\chi(x)-\chi(z),\xi>} a(\chi(x), \chi(z), \xi) u(z) |J| dz d\xi,$$

where $J$ is the Jacobian of the transformation.

Now using Taylor's theorem,

$$\chi(x) - \chi(z) = A(x,z)(x-z)$$

where $A(x,z)$ is a square matrix varying smoothly with $x$ and $z$ and most importantly

$$A(x,x) = \frac{\partial \chi(x)}{\partial x}.$$

This means that $A$ is invertible on the diagonal and also therefore nearby. (We need only worry about behaviour close to the diagonal as the phase is clearly critical only there and so the kernel is smooth off the diagonal.)

So our phase is

$$< \chi(x) - \chi(z), \xi > = < x - z, A^t(x,z)\xi > .$$

So if we do another change of variables, $\eta = A^t(x,z)\xi$,

$$P_1 u(x) = \int e^{i<x-z,\eta>} a(\chi(x), \chi(z), (A^t)^{-1}\eta) u(z) |K| |J| dz d\eta, \quad (8)$$

where $K$ is the determinant of this transformation. This is roughly in the form of a pseudo-differential operator. Both $|K|, |J|$ are smooth functions independent of $\eta$ and so will be harmless. Our problem is thus to check that

$$a_1(x, z, \eta) = a(\chi(x), \chi(z), (A^t)^{-1}(x,z)\eta)$$

is a symbol.

To prove this we use a lemma which often comes in handy when trying to establish something is a symbol.

**Lemma 8.1.** *If $\mathcal{F}_m \subset C^\infty(\mathbb{R}^n_x \times \mathbb{R}^k_\theta)$ for all $m$ and if $b \in \mathcal{F}_m$ implies that*

$$\frac{\partial b}{\partial x_j} = \sum_\alpha a_\alpha b_\alpha + \sum_\alpha a'_\alpha b'_\alpha$$

*with $a_\alpha \in S^0, a'_\alpha \in S^1$ and $b_\alpha \in \mathcal{F}_m, b'_\alpha \in \mathcal{F}_{m-1}$ and*

$$\frac{\partial b}{\partial \theta_k} = \sum_\beta c_\beta d_\beta$$

*with $c_\beta \in S^0$ and $d_\beta \in \mathcal{F}_{m-1}$ and for each element, $b$, of $\mathcal{F}_m$ and compact $K$ there exists $C_K$ such that*

$$|b(x, \theta)| \leq C_K < \theta >^m, \; x \in K, K \; compact$$

*then*

$$\mathcal{F}_m \subset S^m(\mathbb{R}^n; \mathbb{R}^k).$$

*Proof.* This follows easily from repeated use of the Leibniz rule for products. $\square$

So now defining $\mathcal{F}_m = \{a(\chi(x), \chi(z), (A^t)^{-1}(x, z)\eta) : a \in S^m\}$, it follows that $a_1$ is a symbol of order $m$ and so $P_1$ is a pseudo-differential operator of order $m$.

As usual, what we really want to know is how the principal symbol transforms and we also want to show that $P_1$ is classical if $P$ is.

So if we take $a$ to be independent of $y$, which we can do, then the total symbol of $P_1$ is

$$a(\chi(x), (A^t)^{-1}(x, z)\eta)|J||K|.$$

We can now compute the left symbol of $P_1$ in the usual way:

$$\sigma_L(P_1) \sim \sum_\alpha \frac{1}{\alpha!} \partial_z^\alpha \left( D_\eta^\alpha a(\chi(x), (A^t)^{-1}(x, z)\eta)|J||K| \right)_{z=x}.$$

It is now obvious that $\sigma_L(P_1)$ is classical and we have a way of computing the total symbol. If $a_m$ is the top order part of $a$ then the principal symbol of $P_1$ will be

$$a_m(\chi(x), (A^t)^{-1}(x, x)\eta)|J(x, x)||K(x, x)|.$$

However $|J(x, x)||K(x, x)| = 1$ as these are the Jacobians of opposite transformations and

$$(A^t)^{-1}(x, x) = \left( \frac{\partial \chi}{\partial x}^t \right)^{-1}.$$

This says precisely that the principal symbol of $P_1$ is the principal symbol of $P$ pulled back by the lifting of $\chi$ to a function on the cotangent bundle.

We can use this to deduce the coordinate invariance of $H_c^s$. If $P_s$ is a classical, properly-supported, elliptic, pseudo-differential operator of order $s$ then $u$ in $H_c^s$ if and only if $Pu \in L_c^2$. But we know $L_c^2$ is coordinate invariant and $P_s$ will transform to an operator with the same properties so coordinate invariance is immediate.

## 9. Wavefront Sets

Recall the definition of wavefront set

**Definition 9.1.** *If $u \in \mathcal{D}'(\mathbb{R}^n)$ then*

$$(x_0, \xi_0) \notin \mathrm{WF}(u) \subset \mathbb{R}^n \times (\mathbb{R}^n - \{0\})$$

*if and only if there exists $C' > 0$, and $\phi$ smooth and of compact support such that $\phi(x_0) \neq 0$ and for all $N$*

$$|\widehat{\phi u}(\xi)| \leq C_N < \xi >^{-N}$$

*for*

$$\left| \frac{\xi}{|\xi|} - \frac{\xi_0}{|\xi_0|} \right| \leq C'.$$

We can reexpress this using pseudo-differential operators but first we need to introduce the concept of the wavefront set of a pseudo-differential operator! (note this is not the same as the wavefront set of its Schwarz kernel although they are closely related.)

**Definition 9.2.** *If $P \in \Psi DO^m$ then $(x_0, \xi_0) \notin \mathrm{WF}'(P)$ if and only if $P$ has total left symbol $p(x, \xi)$ such that*

$$|D_x^\alpha D_\xi^\beta p(x, \xi)| \leq C_{N, \alpha, \beta} < \xi >^{-N}$$

*for $x$ near $x_0$ and $\xi$ in a conic neighbourhood of $\xi_0$.*

*Remark* 9.3. In fact the estimate for $\alpha = \beta = 0$ together with the fact that $p$ is a symbol is enough. Try to prove this!

The idea here is that $\mathrm{WF}'(P)^c$ expresses the directions in which $P$ kills all singularities as it is of order $-\infty$. Note that using the right symbol gives the same wavefront set. We trivially have that

$$\mathrm{WF}'(AB) \subset \mathrm{WF}'(A) \cap \mathrm{WF}'(B).$$

-once everything has been killed in a certain direction, we will not be able to recover it.

Note also that $P \in \Psi DO^m$ will be $\Psi DO^{-\infty}$ if and only if $\mathrm{WF}'(P) = \emptyset$.

This leads to the important concept of micro-ellipticity which we will show is equivalent to inverting up to order $-\infty$ in a given direction.

**Definition 9.4.** *If $P \in \Psi DO_{cl}^m$ then $P$ is micro-elliptic of order $m$ at $(x_0, \xi_0)$ if $\sigma_m(P)(x_0, \xi_0) \neq 0$.*

**Theorem 9.5.** *If $P \in \Psi DO_{cl}^m$ then $P$ is micro-elliptic at $(x_0, \xi_0)$ if and only if there exists $Q \in \Psi DO_{cl}^{-m}$ such that*

$$(x_0, \xi_0) \notin \mathrm{WF}'(PQ - \mathrm{Id}), (x_0, \xi_0) \notin \mathrm{WF}'(QP - \mathrm{Id}).$$

*$Q$ is said to be a micro-parametrix for $P$.*

*Proof.* Just repeat the proof of the construction of a parametrix of a elliptic operator. The symbols involved may become singular away from $(x_0, \xi_0)$ but if we multiply them all by

$$\phi(x - x_0)\phi\left(\frac{\xi}{|\xi|} - \frac{\xi_0}{|\xi_0|}\right)$$

where $\phi$ is a bump function supported sufficiently close to zero, they will be smooth and retain the property of invertibility close to $(x_0, \xi_0)$. This proves the result. The converse is obvious from considering principal symbols. □

So in some sense the only difficult (interesting) points are where $\sigma_m(P)$ is zero.

**Definition 9.6.** *If $P \in \Psi DO_{cl}^m$ then the characteristic variety of $P$ is*

$$\mathrm{char}(P) = \sigma_m(P)^{-1}(0).$$

**Theorem 9.7.** *If $u \in \mathcal{D}'(\mathbb{R}^n)$ then*

$$\mathrm{WF}(u) = \cap\{\mathrm{char}\, P : Pu \in C^\infty\}.$$

*Proof.* If $(x_0, \xi_0) \notin \mathrm{WF}(u)$ there exists $\phi$ such $\hat{\phi u}$ is rapidly decreasing in direction $\xi_0$ so picking $\psi \in C^\infty(S^{n-1})$ such that

$$\psi\left(\frac{\xi_0}{|\xi_0|}\right) \neq 0$$

and $\hat{\phi u}$ is rapidly decreasing in all directions in supp($\psi$), we pick a pseudo-differential operator $P$ with total symbol

$$\phi(y)\psi\left(\frac{\xi}{|\xi|}\right).$$

It is then clear that

$$Pu \in C^\infty \text{ and } (x_0, \xi_0) \notin \text{WF}(Pu).$$

So we have proven the result one way.

So now suppose $P$ is micro-elliptic at $(x_0, \xi_0)$ and $Pu \in C^\infty$ then applying a micro-parametrix for $P$ at $(x_0, \xi_0)$, we have

$$u + Ru \in C^\infty$$

and $(x_0, \xi_0) \notin \text{WF}'(R)$. So pick $\phi, \psi$ to be bump functions such that

$$\text{supp}\left(\phi(x - x_0)\psi\left(\frac{\xi}{|\xi|} - \frac{\xi_0}{|\xi_0|}\right)\right) \cap \text{WF}'(R) = \emptyset.$$

Then picking $S$ to have total symbol $\phi(y - x_0)\psi\left(\frac{\xi}{|\xi|} - \frac{\xi_0}{|\xi_0|}\right)$, we have that $\text{WF}'(S) \cap \text{WF}'(R) = \emptyset$ and so

$$Su \in C^\infty,$$

(using the pseudo-locality of $S$.)

But this says precisely that

$$(x_0, \xi_0) \notin \text{WF}(u).$$

$\square$

Since we have already shown that the principal symbols of pseudo-differential operators transform as functions on the cotangent bundle, we have that wavefront set is coordinate invariant if we regard it as a subset of the cotangent bundle.

**Theorem 9.8.** *If $P \in \Psi DO_{cl}^m$ is properly supported and $u \in \mathcal{D}'(\mathbb{R}^n)$ then*

$$\text{WF}(Pu) \subset \text{WF}(u) \cap \text{WF}'(P).$$

This says that pseudo-differential operators are *micro-local.* This is a stronger property than pseudo-locality but of course weaker than locality.

*Proof.* To prove the first inclusion, if $(x_0, \xi_0) \notin \mathrm{WF}(u)$ then there exists $Q'$ micro-elliptic at $(x_0, \xi_0)$ such that $Q'u$ is smooth. Taking a micro-local parametrix, $Q''$, for $Q'$ at $(x_0, \xi_0)$, we have putting $Q = Q''Q'$ that $Qu \in C^\infty$ and that the *total* symbol of $Q$ is identically one near $(x_0, \xi_0)$ i.e. $Q = \mathrm{Id} + R$ with $(x_0, \xi_0) \notin \mathrm{WF}'(R)$. Now, $PQu$ will of course be smooth and

$$QPu = PQu + [Q, P]u.$$

But $(x_0, \xi_0) \notin \mathrm{WF}'([Q, P])$ as

$$[Q, P] = (\mathrm{Id} + R)P - P(\mathrm{Id} + R)$$
$$= RP - PR.$$

So picking $S$ micro-elliptic at $(x_0, \xi_0)$ such that

$$\mathrm{WF}'(S) \cap \mathrm{WF}'(R) = \emptyset$$

we have that

$$SQPu \in C^\infty.$$

As $SQ$ is micro-elliptic at $(x_0, \xi_0)$, the first part of the result follows.

The second part is easier. If

$$(x_0, \xi_0) \notin \mathrm{WF}'(P)$$

then we can take $S$ micro-elliptic at $(x_0, \xi_0)$ such that

$$\mathrm{WF}'(SP) \subset \mathrm{WF}'(P) \cap \mathrm{WF}'(S) = \emptyset$$

and so $SP \in \Psi DO^{-\infty}$ This implies that

$$SPu \in C^\infty$$

which proves the result. $\qquad\qquad\qquad\qquad\qquad\qquad\qquad\square$

**Corollary 9.9.** *If $P \in \Psi DO_{cl}^m$, properly supported then*

$$\mathrm{WF}(u) \subset \mathrm{WF}(Pu) \cup \mathrm{char}(P).$$

*Proof.* This is just a simple application of micro-ellipticity. If

$$(x_0, \xi_0) \notin \mathrm{WF}(Pu) \cup \mathrm{char}(P)$$

then $P$ is micro-elliptic at $(x_0, \xi_0)$ and $(x_0, \xi_0) \notin \mathrm{WF}(Pu)$. So letting $Q$ be a micro-local parametrix for $P$ at $(x_0, \xi_0)$, we have that

$$(x_0, \xi_0) \notin \mathrm{WF}(QPu).$$

But

$$QPu = u + Ru$$

with

$$(x_0, \xi_0) \notin \mathrm{WF}'(R)$$

and hence

$$(x_0, \xi_0) \notin \mathrm{WF}(Ru).$$

So

$$(x_0, \xi_0) \notin \mathrm{WF}(u)$$

which is what we wanted to prove.                                         $\square$

## 10. The Propagation of Singularities

In the previous section, we showed that if $Pu \in C^\infty$ then

$$\mathrm{WF}(u) \subset \{\sigma_m(P) = 0\}.$$

The obvious next question is which subsets are possible? We look at this question in the case when $\sigma_m(P)$ is real. The answer due to Hörmander is that the wavefront set is a union of integral curves of the Hamiltonian of the principal symbol of $P$. If $p_m = \sigma_m(P)$ then the Hamiltonian is the vector field on $T^*(\mathbb{R}^n)$

$$H_{p_m} = \sum_{j=1}^{n} \frac{\partial p_m}{\partial \xi_j} \frac{\partial}{\partial x_j} - \frac{\partial p_m}{\partial x_j} \frac{\partial}{\partial \xi_j}$$

or if you prefer

$$\left( \frac{\partial p_m}{\partial \xi_1}, \frac{\partial p_m}{\partial \xi_2}, \dots, \frac{\partial p_m}{\partial \xi_n}, -\frac{\partial p_m}{\partial x_1}, \dots, -\frac{\partial p_m}{\partial x_n} \right).$$

*Remark* 10.1. The first order differential operator and the vector field are usually identified in the study of differential geometry. The point being that if $\gamma(t)$ is an integral curve of a vector field $v = (v_1, \dots, v_n)$ on $\mathbb{R}^n$ then

$$\frac{d}{dt}(f \circ \gamma)(t) = \sum v_i(\gamma(t)) \frac{\partial f}{\partial x_i}(\gamma(t)).$$

A good question is, of course, whether this is the end of the story. For example, if the wavefront set had to always be the entire characteristic variety or empty then the result would be true but boring. In fact, Hörmander also showed that any closed subset of the characteristic variety which is a union of integral curves is possible but that is beyond the scope of this course.

**Theorem 10.2.** *Let $P$ be a proper, classical pseudo-differential operator and let $u \in \mathcal{D}'$ be such that $Pu \in C^\infty$ then $\mathrm{WF}(u)$ is a union of integral curves of the Hamiltonian of $\sigma_m(P)$.*

There are two basic approaches to proving this result the highbrow way is to develop the calculus of Fourier integral operators and then reduce to the case where $P = D_{x_1}$ in which case it is easy. The lowbrow and original proof is the method of commutators which we will use here. The main disadvantage of this is that we will need an excursion into functional analysis.

Before we proceed note that

$$H_{p_m} p_m = 0$$

that is $p_m$ is constant along integral curves of $H_{p_m}$, in particular the characteristic variety of $P$ is invariant under the flow. We call the integral curves in the characteristic variety bicharacteristics. Note that this also shows the characteristic variety is a union of bicharacteristics.

If $Q$ is an elliptic differential operator of order $1 - m$ then

$$\mathrm{WF}(QPu) = \mathrm{WF}(u)$$

and if $Q$ has principal symbol $q$ then

$$H_{qp_m} = qH_{p_m} + p_m H_q$$

so in $p_m = 0$, we have that $H_{qp_m} = qH_{p_m}$ that is that $H_{qp_m}$ is a non-zero multiple of $H_{p_m}$ and thus will have the same integral curves. It is thus enough to prove the theorem in the special case when $P$ is of first order.

As the WF set is a closed set it is enough, by connectedness, to show that if a point is in the WF set then it is nearby on the bicharacteristic too so our although our result is global it can be proved locally. We thus fix some point $x_0$ and let $\phi$ be compactly supported and identically one near $x_0$ and $\psi$ be identically one near $x_0$ supported within the set where $\phi$ is identically one. Then if $Pu$ is smooth so is $\psi P(\phi u)$. We then have that $\psi P \phi$ has compact kernel and that it has the same bicharacteristics in a neighbourhood of $x_0$ - this means it is enough to consider the case where $P$ has a compact kernel.

In general two pseudo-differential operators $P \in \Psi DO^m_{cl}, Q \in \Psi DO^{m'}_{cl}$ will not commute eg

$$P = D_{x_1}, Q = x_1$$

but the principal symbols

$$\sigma_{m+m'}(PQ) = \sigma_m(P)\sigma_{m'}(Q) = \sigma_{m+m'}(QP)$$

are equal so

$$[P, Q] = PQ - QP \in \Psi DO^{m+m'-1}$$

and consulting our asymptotic expansions has principal symbol $\frac{1}{i}\{p_m, q_{m'}\}$ that is

$$\frac{1}{i}\left(\sum_{j=1}^{n}\frac{\partial p_m}{\partial \xi_j}\frac{\partial q_{m'}}{\partial x_j} - \frac{\partial p_m}{\partial x_j}\frac{\partial q_{m'}}{\partial \xi_j}\right) = \frac{1}{i}H_{p_m}q_{m'}$$

which is the crucial point here.

We will prove the result by considering the operator $D_t - P$ and constructing operators which commute with it. If $(x_0, \xi_0) \notin \mathrm{WF}(u)$ then there is a $Q' \in \Psi DO_{cl}^0(\mathbb{R}^n)$ micro-elliptic at $(x_0, \xi_0)$ such that $Q'u \in C^\infty$. We try to find a family $Q(t, x, D_x)$ which are pseudo-differential for each $t$, with total symbol varying smoothly with $t$ such that

$$(D_t - P)Q = Q(D_t - P)$$

and $Q'(x, D_x) = Q(0, x, D_x)$. It then follows that

$$(D_t - P)(Qu) \in C^\infty$$

and

$$Qu_{|t=0} \in C^\infty.$$

There is a regularity theorem which we will prove that then guarantees that $Qu$ is smooth for all $t$. This means that

$$\mathrm{WF}(u) \subset \{q_0(t, x, \xi) = 0\}$$

for each $t$. Let $\gamma_{x_0, \xi_0}(t)$ be the bicharacteristic through $(x_0, \xi_0)$ then we will show that $q_0(t, \gamma_{(x_0, \xi_0)}(t))$ is zero if and only if $q_0(0, x_0, \xi_0)$ is zero and the theorem then follows.

So how do we construct $Q$ ? The principal symbol of $[D_t - P, Q]$ is

$$\frac{\partial q_0}{\partial t} - H_{p_1}q_0.$$

We want this to be zero and $q_0(0, x, \xi) = q(x, \xi)$. This is just a standard first order linear PDE and is solved by making $q$ constant along the integral curves of $(1, H_{p_1})$ which are of the form $(t, \gamma(-t))$ where $\gamma(t)$ is an integral curve of $H_{p_1}$. So $q$ is either 0 or non-zero all along such an integral curve which is what we need.

We also need to construct the lower order terms of $Q$. If we have constructed the higher terms then the next term in the symbol expansion will be of the form

$$\frac{\partial q_{-j}}{\partial t} - H_{p_1}q_{-j} = \text{ function of higher order terms}$$

this being linear we can as usual solve. Asymptotically summing we then have $Q$ which commutes with $D_t - P$ up to smooth terms which is enough.

So we have proven the result modulo the proof of the regularity theorem. We prove the regularity result by using Sobolev spaces. Since we have restricted to operators $P$ which have compact kernels and real principal symbol we have that $P - P^*$ is zeroth order and has a compact kernel - this guarantees that it is continuous on all Sobolev spaces.

We first prove an energy estimate and then use a duality argument to construct a solution of $D_t - P$ with the right regularity properties which by uniqueness then agrees with the original solution which means it has the right regularity properties too.

**Proposition 10.3.** *Let $L = D_t - P$ then for all $s \in \mathbb{R}$, there exists $\lambda_s$ such that for $\lambda \geq \lambda_s$, we have for*

$$u \in C^1\left([0,T]; H^s\right) \cap C^0\left([0,T], H^{s+1}\right)$$

*that*

$$\sup_{t \in [0,T]} e^{-\lambda t} \|u(t,.)\|_s \leq \|u(0,.)\|_s + 2 \int_0^T e^{-\lambda t} \|Lu(t,.)\|_s dt.$$

*Proof.* It will be slightly more convenient to eliminate the $1/i$ so we instead study

$$L' = \frac{\partial}{\partial t} - Q$$

with $Q = iP$ so $R = Q + Q^*$ is continuous on $H^s$.

We first do the case with $s = 0$ that is for $L^2$. We have

$$\frac{d}{dt} \|e^{-\lambda t} u(t,x)\|_0^2 = \frac{d}{dt} \int e^{-2\lambda t} u(t,x) \bar{u}(t,x) dx$$

which equals

$$-2\lambda \|e^{-\lambda t} u(t,x)\|_0^2 + 2\Re \langle e^{-2\lambda t} \frac{du}{dt}, u \rangle.$$

Now

$$\langle Qu, u \rangle = \langle u, Q^* u \rangle = -\langle u, Qu \rangle + \langle u, Ru \rangle = -\overline{\langle Qu, u \rangle} + \langle u, Ru \rangle,$$

which implies that

$$2\Re \langle Qu, u \rangle = \langle u, Ru \rangle \leq \|R\| \|u\|_0^2$$

where $\|R\|$ is the operator norm on $L^2$. So we conclude that

$$2\Re \langle e^{-2\lambda t} L'u, u \rangle \geq \frac{d}{dt} \left\| e^{-\lambda t} u(t,.) \right\|_0^2 + (2\lambda - \|R\|) \left\| e^{-\lambda t} u(t,.) \right\|_0^2.$$

So provided $\lambda \geq \|R\|/2$, we can drop the last term and then integrating from 0 to $t$ we have that

$$\left\| e^{-\lambda t} u(t, .) \right\|_0^2 - \left\| e^{-\lambda t} u(0, .) \right\|_0^2 \leq 2 \int\limits_0^T \Re \langle e^{-2\lambda t} L'u, u \rangle dt.$$

Putting $M = \sup\limits_{t \in [0,T]} \left| e^{-\lambda t} u(t, .) \right|_0$ and applying Cauchy Schwartz to $\Re \langle e^{-2\lambda t} L'u, u \rangle$, we have that

$$M^2 \leq \|u(0, .)\|^2 + 2M \int\limits_0^T e^{-\lambda t} \|Lu(t, .)\|_0 dt.$$

As $\|u(0, .)\|_0 \leq M$, we have, dividing by $M$, that the result follows for $L^2$.

To do the general case, we conjugate by the operators $\Lambda_s$ from section 6. We then have

$$Lu = \Lambda_{-s} \tilde{L} \Lambda_s u$$

and $\tilde{L}$ is obtained from $L$ by replacing $P$ by $\tilde{P} = \Lambda_s P \Lambda_{-s}$. So

$$\|Lu\|_s = \|\tilde{L}\Lambda_s u\|_0.$$

<div style="text-align: right">□</div>

**Theorem 10.4.** *Let $L$ be as above and $s \in \mathbb{R}$. For every $f \in L^1((0, T); H^s)$ and $\phi \in H^s$ there is then a unique solution $u \in C([0, T]; H^s)$ of*

$$Lu = f, u_{|t=0} = \phi$$

*and this solution satisfies the energy estimate* (10.3).

*Proof.* We work with $L' = iL$ as before.

To see uniqueness that if $u_1, u_2$ are both solutions then $u = u_1 - u_2$ satisfies the equation with $f, \phi = 0$. We have that

$$\frac{\partial u}{\partial t} = Pu \in C([0, T]; H^{s-1})$$

and so we have that

$$u \in C^1([0, T]; H^{s-1}) \cap C^0([0, T]; H^s).$$

So applying Proposition 10.3 with $s$ replaced by $s - 1$ we have that

$$\sup e^{-\lambda t} \|u(t, .)\|_{(s-1)} = 0$$

and thus that $u$ is zero.

To construct a solution we use a duality argument and the Hahn-Banach theorem. If $u$ is a solution of the Cauchy problem and $v \in C_0^\infty((t,x) : t < T)$ then

$$\int_0^T \langle u, -\frac{\partial v}{\partial t} + P^* v \rangle dt = \int_0^T \langle f, v \rangle dt + \langle \phi(.), v(0,.) \rangle.$$

The idea here is to extend the functional defined by the left hand side on elements of the form $L^* v$ to the entire space of functions $L^1 \left( (0,T); H^{-s} \right)$ by Hahn-Banach to obtain an element $u$ of

$$\left( L^1 \left( (0,T); H^{-s} \right) \right)^* = L^\infty \left( [0,T]; H^s \right)$$

which will then solve the problem.

To apply the Hahn-Banach theorem we will need to show that the functional is well-defined and continuous. Now applying Prop 10.3 to $\tilde{v}(t) = v(T' - t)$, with 's = -s' for any $T' < T$ we have that

$$\sup_{t \in [0,T']} \|e^{-\lambda t} \tilde{v}(t,.)\|_{-s} \leq C \int_0^{T'} e^{-\lambda t} \|g\|_{-s} dt$$

with

$$g = \frac{\partial \tilde{v}}{\partial t} + P \tilde{v}.$$

So for any $T'$, $\|v(T - T')\|_{-s}$ is dominated by the right hand side so

$$\sup_{t \in [0,T]} \|v\|_{-s} \leq C \int_0^T \|g\|_{-s} dt.$$

This shows that $L^*$ is injective so we can define a map from $L^*(C_0^\infty((t,x) : t < T))$ to $\mathbb{C}$ by

$$\tilde{\psi}(L^* v) = \int_0^T \langle f(t,.), v(t,.) \rangle dt + \langle \phi, v(0,.) \rangle.$$

The estimate also shows that

$$\|\tilde{\psi}(L^* v)\| \leq C \left( \|\phi\|_s + \int_0^T \|f(t,.)\|_s dt \right) \int_0^T \|L^* v(t,.)\|_{-s} dt$$

that is that $\tilde{\psi}$ is continuous on $L^*(C_0^\infty((t,x) : t < T))$ if one regards it as a subspace of $\left( L^1((0,T), H^{-s}) \right)$.

The Hahn-Banach says that one can extend $\tilde{\psi}$ to a continuous linear functional on $\left( L^1((0,T), H^{-s}) \right)$ that is to an element of

$$L^\infty((0,T), H^s).$$

Call this element $u$ and we have

$$\int \langle u, -\frac{\partial v}{\partial t} + Q^* v \rangle dt = \int\limits_0^T \langle f(t,.), v(t,.) \rangle dt + \langle \phi, v(0,.) \rangle$$

for $v \in C_0^\infty((t,x) : t < T)$. So certainly, on $(0,T)$ we have

$$\frac{\partial u}{\partial t} - Qu = f$$

in a distributional sense. Now if we assume for now that $f$ is Schwartz then we have that

$$\frac{\partial u}{\partial t} = Qu + f \in L^\infty((0,T); H^{s-1}).$$

So $u$ is in $C([0,T]; H^{s-1})$ and then if we assume $\phi$ is Schwartz we have integrating the equation again that

$$u \in C^1([0,T]; H^{s-2})$$

and has the right initial value - $\phi$.

For $f, \phi$ Schwartz the same argument would work with $s$ replaced by $s + 2$ so in that case we have a solution $u'$ which will equal $u$ by uniqueness in $C^1([0,T]; H^s)$ and will satisfy the energy estimate (10.3).

To do the general case, we use the density of Schwartz space in $H^s$ and $L^1((0,T); H^s)$. (Prove this!) Let $f_k \to f$ in $L^1((0,T); H^s)$, $f_k$ Schwartz and $\phi_k \to \phi$, $\phi_k$ Schwartz. Then for each $k$ we can use the previously proven to case to get

$$u_k \in C^1([0,T]; H^s) \cap C^0([0,T]; H^{s+1})$$

solving with data $(f_k, \phi_k)$. But then applying (10.3) again, we deduce that $u_k$ is Cauchy and that the limit is the desired solution $u$. $\qquad \square$

So why does all this allow us to prove the propagation of singularities? We have reduced to the case where $u \in \mathcal{E}'(\mathbb{R}^n)$ and $P$ has compactly supported kernel. We take $Q(t,x,D_x)$ commuting up to smoothing with $D_t - P$ as above and $Q(0,x,D_x)u \in C_0^\infty$. We also assume that $Q$ is uniformly properly supported.

As $u$ is compactly supported, it must be in some Sobolev space - just estimate

$$\langle u, e^{-ix.\xi} \rangle.$$

Since we have chosen $Q$ to be a smooth family of pseudo-differential operators, we have

$$Qu \in C^\infty([0,T]; H^s)$$

and

$$(D_t - P)Qu \in C_0^\infty.$$

If we take $s' > s$ then the theorem says there is $v$ satisfying the equation with the same initial value and $v(t,.) \in H^{s'}$. But by uniqueness the two solutions agree so $Qu(t,.)$ is in $H^{s'}$. But $s'$ was arbitrary so we conclude that $Qu(t,.)$ is smooth for each $t$ which proves our result.

## 11. Operators on Manifolds

So having developed a theory of pseudo-differential operators on $\mathbb{R}^n$ which is coordinate invariant, we can reduce the manifold case to $\mathbb{R}^n$ by working locally.

Let $X$ be a smooth manifold with an atlas $\{U_\alpha, V_\alpha, \phi_\alpha\}$.

**Definition 11.1.** *The map $P$ from $C_0^\infty(X)$ to $C^\infty(X)$ is a (classical) pseudo-differential operator of order $m$ if in any local coordinate patch $U_\alpha$, with $\psi, \phi \in C_0^\infty(U_\alpha)$ $\psi P \phi$ defines a (classical) pseudo-differential operator in local coordinates $C_0^\infty(U_\alpha) \to C_0^\infty(U_\alpha)$ of order $m$. We also require that if $\phi, \psi$ are in $C_0^\infty(X)$ with disjoint supports then $\phi P \psi$ has a smooth Schwartz kernel.*

The second condition is a decree that our operator be pseudo-local. Why do we need this? Well coordinate patches are local objects so they say nothing about behaviour between points far away. For example if one takes a sphere and define the map

$$Af(x) = f(-x)$$

then if one takes an atlas with small coordinate patches, it is zero locally and hence obeys the first condition. But it is not a pseudo-differential operator.

Since we previously proved coordinate invariance, if an operator is pseudo-differential with respect to some atlas then it is pseudo-differential with respect to any atlas.

We now want to extend our symbol calculus from $\mathbb{R}^n$ to $X$. The total symbol is a nasty object - evidenced by our formulas for the change of coordinates but the principal symbol is nice. Let $SH^m(T^*(X) - \{0\})$ be smooth functions on $T^*(X) - \{0\}$ which are positively homogeneous of degree $m$. (Recall that the cotangent bundle is a vector bundle so on each fibre there is a natural scalar multiplication.) To define the principal symbol near a point $x$ we take a coordinate chart $(U_\alpha, V_\alpha, \phi_\alpha)$ and a bump function $\psi$ identically one near $x$ and define $\sigma_m(P)(x, \xi)$ to be the principal symbol of $\phi P \phi$ in the local coordinates on the cotangent bundle naturally induced by the local coordinates on the manifold. We then have

**Proposition 11.2.** *There is a natural short exact sequence,*

$$0 \to \Psi DO_{cl}^{m-1}(X) \to \Psi DO_{cl}^m(X) \to SH^m(T^*(X) - \{0\}) \to 0.$$

The rest of our results we then take over from the $\mathbb{R}^n$ case.

**Proposition 11.3.** *If $P, Q$ are properly supported pseudo-differential operators of order $m, m'$ then $PQ$ is a pseudo-differential operator of order $m + m'$ and*

$$\sigma_{m+m'}(PQ) = \sigma_m(P)\sigma_{m'}(Q).$$

The asymptotic completeness of the calculus of pseudo-differential operators will also hold.

**Proposition 11.4.** *If $P_j \in \Psi DO_{cl}^{m-j}(X)$ for $j = 0 \ldots \infty$ then there exists $P \in \Psi DO_{cl}^m(X)$ such that*

$$P - \sum_{j < N} P_j \in \Psi DO_{cl}^{m-N}(X)$$

*for all $N$.*

*Proof.* Let $f_\alpha$ be a partition of unity subordinate to some atlas $(U_\alpha, V_\alpha, \phi_\alpha)$ and let $g_\alpha \in C_0^\infty(U_\alpha)$ be identically one on the support of $f_\alpha$. We then have in local coordinates that

$$P_{\alpha,j} = g_\alpha P_j f_\alpha \in \Psi DO_{cl}^{m-j}(\mathbb{R}^n)$$

and we can asymptotically sum to get $P_\alpha$. We then just put

$$P = \sum g_\alpha P_\alpha g_\alpha$$

and we are done. We use the cut-offs to ensure that the operators have compact kernels contained in $U_\alpha \times U_\alpha$. (we can assume the cover $U_\alpha$ is locally finite.) $\qquad \square$

We define ellipticity as before - the principal symbol does not vanish. We can then construct a parametrix for an elliptic operator - there are two basic approaches. The first is work locally and patch together - this is fiddly and unnecessarily so. The second is to use the principal symbol map to get a $Q \in \Psi DO_{cl}^{-m}$ with

$$\sigma_{-m}(Q) = \sigma_m(P)^{-1}$$

so

$$PQ = I - R, \ R \in \Psi DO_{cl}^{-1}(X)$$

and then let

$$Q' \sim \sum R^j.$$

This will then be a parametrix. If one needed to, one could compute the total symbol in any local coordinates in the same way as for $\mathbb{R}^n$.

We can also define wavefront set as before.

**Definition 11.5.** *If $u \in \mathcal{D}'(X)$ then*

$$\mathrm{WF}(u) = \bigcap_{Pu \in C^\infty} \mathrm{char}(P) \subset T^*(X) - \{0\}.$$

We then have micro-locality, micro-ellipticity and so on. Propagation of singularities will also hold as its make invariant sense on the cotangent bundle - our proof was local anyway.

The continuity on $L^2_c(X)$ of properly supported operators will go over since one can cover any compact set by a finite number of coordinate patches and use a partition of unity.

So the theory for a general manifold is not very different from the theory we have developed for $\mathbb{R}^n$. However if we work on a compact manifold, we can get some extra results - we can show that elliptic operators are Fredholm either as operators on smooth functions or between Sobolev spaces. Note that elliptic regularity guarantees that the kernel consists of smooth functions and thus that it will be independent of which Sobolev space we work on.

We want to define Sobolev spaces on a compact manifold. We showed previously that an element is in $H^s_c(\mathbb{R}^n)$ if and only if $P_s u \in L^2_c(\mathbb{R}^n)$ for $P_s$ elliptic of order $s$ and properly supported. We also know that $C^\infty_0(\mathbb{R}^n)$ is dense in $H^s_c$ so we define for $X$ a compact manifold, and $s > 0$, that $H^s(X)$ is the completion of $C^\infty(X)$ with respect to the inner product

$$\langle u, v \rangle_s = \langle u, v \rangle_{L^2} + \langle P_s u, P_s v \rangle_{L^2},$$

where $P_s$ is a self-adjoint elliptic pseudo-differential operator of order $s$ on $X$. We define for $s < 0$ by duality.

The first term is required in case $P_s$ has non-trivial kernel (consisting of smooth functions.) This definition makes $H^s$ a Hilbert space with norm dependent on $P_s$ but a different choice will give an equivalent norm. Suppose we took $Q_s$ instead then $Q_s = R_0 P_s + E$ with $E$ smoothing and $R_0$ zeroth order and thus continuous on $L^2$. We then have that the norm squared coming from $Q_s$ is

$$\langle u, u \rangle_{L^2} + \langle R_0 P_s u, R_0 P_s u \rangle_{L^2} + 2\Re\langle Eu, R_0 P_s u \rangle_{L^2} + \langle Eu, Eu \rangle_{L^2}$$

using the continuity of $E, R_0$ on $L^2$ and Cauchy-Schwartz this is less than or equal to

$$C(\|u\|_{L^2} + \|P_s u\|)^2.$$

By symmetry, we have that the norms are equivalent and thus a natural topology.

Of course, we want that if $P$ is order $m$ then

$$P : H^s \to H^{s-m}.$$

For zeroth order operators on $L^2$, this is trivial as above. We then write

$$P = P_{m-s}(P_{s-m} P P_{-s}) P_s + R$$

with $R$ order $-\infty$. Provided $s \geq m$ this then immediately gives the continuity. We deduce the case on negative Sobolev spaces by duality. To stitch the two together, we can write $R$ as a product of operators, (up to smoothing at least) one mapping into $L^2$ and the other mapping out of it.

If $u$ is supported in a coordinate patch, we have two ways to compute its $H^s$ norm - its Sobolev norm in local coordinates and a norm on the manifold. As the $L^2$ norms are equivalent and $P_s$ is continuous, with continuous parametrix, on both sorts of Sobolev spaces, we see that these are equivalent. This gives us an alternate definition of the norm on the manifold, let $\sum \phi_\alpha$ be a finite partition of unity with $\phi_\alpha$ supported in a coordinate patch then let

$$\|u\|_s = \sum_\alpha \|\phi_\alpha u\|_s$$

where $\|\phi_\alpha u\|_s$ is taken either in local coordinates or on the manifold as they are equivalent anyway. This is equivalent to the original norm - multiplication by $\phi_\alpha$ is continuous as a zeroth order pseudo-differential operator so $\|\phi_\alpha u\| \leq C\|u\|$ and thus summing gives the inequality one way; the other way comes the triangle inequality applied to $\sum \phi_\alpha u$.

The important thing about Sobolev spaces on a compact manifold or in a compact set on $\mathbb{R}^n$ is Rellich's compactness lemma:

**Theorem 11.6.** *If $t > s$ and $X$ is a compact manifold then the inclusion*

$$i : H^t(X) \to H^s(X)$$

*is compact.*

*Proof.* We first reduce to the compact set in $\mathbb{R}^n$ case. We need to show a bounded sequence in $H^t(X)$ will have a convergent subsequence in $H^s(X)$.

Let $\phi_i, i = 1 \ldots k$ be a partition of unity subordinate to an atlas and let $u_n$ be a bounded sequence in $H^t(X)$. We then have that

$$\phi_1 u_n$$

has a convergent subsequence from the $\mathbb{R}^n$ case, $\phi_1 u_n^1$. But then $\phi_2 u_n^1$ has one too and so on. The $k^{th}$ subsequence $u_n^k$ will be such that $\phi_j u_n^k$ is convergent for all $k$ and thus is convergent from our result on the equivalence of norms.

So we need only prove the result in $\mathbb{R}^n$. Let $u_n$ be a bounded sequence supported in a compact set $K \subset \mathbb{R}^n$. Now if $\phi$ is identically one on $K$ then $\phi u = u$ for all $u$ supported in $K$. So

$$\hat{u} = \left(\frac{1}{2\pi}\right)^n \hat{u} * \hat{\phi}$$

and thus that

$$D_\xi^\alpha \hat{u} = \left(\frac{1}{2\pi}\right)^n \hat{u} * D_\xi^\alpha \hat{\phi}.$$

As

$$\langle \xi \rangle^s \leq C \langle \xi - \eta \rangle^{|s|} \langle \eta \rangle^s,$$

we deduce that

$$|\langle \xi \rangle^s D_\xi^\alpha \hat{u}(\xi)| \leq \int \langle \eta \rangle^s |\hat{u}(\eta)| \langle \xi - \eta \rangle^{|s|} |D^\alpha \hat{\phi}(\xi - \eta)| d\eta$$

and thus, using Cauchy-Schwartz, is bounded by $C_\alpha \|u\|_s$.

Applying this to the sequence $u_n$ the first partial derivatives are bounded uniformly, independent of $n$, on compact sets. Thus we have that $u_n$ is equicontinuous on compact sets and is a bounded sequence on compact sets. So by the Arzela-Ascoli theorem there is a convergent subsequence with respect to the uniform norm of the Fourier transform on any compact set. Using diagonalization, we can get the same sequence for all compact sets - as $\mathbb{R}^n$ is a countable union of compact sets.

Let $v_k$ be the resulting subsequence. We show that this is Cauchy in $H^s$ and the result then follows.

We can write

$$\|v_k - v_l\|_s^2 = I_1 + I_2$$

with

$$I_1 = \int\limits_{|\xi| \leq R} (1 + |\xi|^2)^s |\hat{v}_k(\xi) - \hat{v}_l(\xi)|^2 d\xi$$

and

$$I_2 = \int\limits_{|\xi| \geq R} (1 + |\xi|^2)^{s-t}(1 + |\xi|)^t |\hat{v}_k(\xi) - \hat{v}_l(\xi)|^2 d\xi$$

which is bounded by

$$(1 + R^2)^{s-t} \|u_k - u_l\|_t^2$$

which will converge to zero as $R \to \infty$ using the originial hypothesis that the sequence is bounded in $H^t$.

So given $\epsilon > 0$, we pick $R$ so that the second integral is less than $\epsilon/2$ and then we can make the first integral arbitrarily small as we know $\hat{v}_k$ is uniformly convergent on $|\xi| \leq R$. $\qquad \square$

**Corollary 11.7.** *If $P$ is a pseudo-differential operator of negative order then $P$ is a compact operator on any Sobolev space.*

So given all this, we now have a real theorem,

**Theorem 11.8.** *Let $P$ be an elliptic operator of order $m$ on a compact manifold, $M$, then $P$ is Fredholm as a map from $H^s(M)$ to $H^{s-m}(M)$. The index of $P$ is independent of $s$ and depends only on the principal symbol of $P$.*

*Proof.* Let $Q$ be a parametrix for $P$ then $PQ - I$ is order $-\infty$ on $H^{s-m}(M)$ and $QP - I$ is of order $-\infty$ on $H^s$. Since operators of negative order are compact, we deduce that $P$ is invertible modulo compact operators and thus is Fredholm. (Atkinson's parametrix criterion.)

The index of $P$ is the dimension of the kernel of $P$ minus the dimension of the kernel of $P^*$. By elliptic regularity, these spaces consist of smooth functions which are thus in every Sobolev space and so will be independent of $s$.

If $P' = P + R$ with $R$ of order $m - 1$ then $R$ is compact as operator from $H^s(M)$ to $H^{s-m}(M)$ as we can decompose into a continuous map into $H^{s-m+1}$ and the compact inclusion from there to $H^{s-m}(M)$. Since a compact perturbation does not affect the index, this yields that $P, P'$ have the same index and thus that the principal symbols determine the index. $\qquad \square$

**Corollary 11.9.** *$P$ has finite dimensional kernel and cokernel as a map on $C^\infty(X)$.*

*Proof.* The kernel on any Sobolev space is precisely the smooth functions such $Pu = 0$ so the kernel of $P$ is finite dimensional. As the space of smooth functions is in any Sobolev space we also have that at most finite dimensional subset of them are not in the image.    □

A self-adjoint elliptic operator will have index zero. We really need to make more precise what we mean by this as there are many concepts of self-adjointness. For a map on $C_0^\infty(\mathbb{R}^n)$ there is a natural adjoint obtained by taking the conjugate of the tranpose of the Schwartz kernel - on a manifold we need to worry about a density term. So we say that an operator is formally self-adjoint if

$$\langle Pu, v \rangle = \langle u, Pv \rangle$$

in the $L^2$ pairing for all $u, v \in C_0^\infty(X)$. Now if $P$ is an elliptic pseudo-differential operator of order $m$ then we know

$$P : H^s(X) \to H^{s-m}(X)$$

is Fredholm if $X$ is a compact manifold. Its index is equal to

$$\dim \ker P - \dim \ker P^*.$$

But $P^*$ is the adjoint as a map from

$$(H^{s-m})^* \to (H^s)^*$$

that is from

$$H^{m-s} \to H^{-s}.$$

The important point is that this map is the extension of $P^*$ the formal adjoint on $C^\infty(X)$ so by elliptic regularity the kernel will be equal to that of the formal adjoint on $C^\infty(X)$ which for formally self-adjoint operator is equal to that of $P$. Hence the index is zero.

Another important property of formally self-adjoint, elliptic operators on compact manifolds, is that they give a natural decomposition of spaces. If we take $v \in C^\infty(X)$ and $w \in \ker P$ (which is smooth so independent of $s$) then

$$\langle Pv, w \rangle_{L^2} = \langle v, Pw \rangle_{L^2} = 0$$

So $P(C^\infty(X))$ is orthogonal in $L^2$ to $\ker P$ but $C^\infty(X)$ is dense in $H^m(X)$ so we have that $P(H^m(X)$ is orthogonal to $\ker P$. However the codimension of $P(H^m(X))$ is equal to the dimension of $\ker P$ so we have that

$$L^2(X) = P(H^m(X)) \oplus \ker P.$$

Invoking elliptic regularity again, the same is true for $C^\infty(X)$. So to summarize

**Proposition 11.10.** *Let $P$ be an elliptic, formally self-adjoint pseudo-differential operator on $X$ then $P$ is of index zero and*

$$C^\infty(X) = P(C^\infty(X)) \oplus \ker P.$$

Note that if two principal symbols were homotopic, we could with a little effort construct a continuous family of operators quantizing them, as index is constant on continuous families we could then deduce that the principal symbols gave the same index and thus that index is really a map from homotopy classes of principal symbols to $\mathbb{Z}$ - ie it is a topological object. In fact, the index is determined by the homotopy class as a map from the cosphere bundle - if we pick a Riemannian metric this is the covectors of modulus one with respect to the dual metric.

**Proposition 11.11.** *Let $X$ be a Riemannian manifold and $m \in \mathbb{R}$. There exists a self-adjoint operator $A_m$ of order $m$ with principal symbol*

$$|\xi|_x^m.$$

*Proof.* Let $B$ be an operator of order $m/2$ with symbol $|\xi|_x^{m/2}$ then $A = B^*B$ has all the requisite properties. (with a smooth perturbation one could actually make $B$ invertible.) $\qquad\square$

Now a self-adjoint operator always has index zero so the index of $P$ is equal to that of $A_{-m}P$. This means that the index of the principal symbol will be equal to that obtained by taking the zeroth order operator with symbol defined by projection on to the cosphere bundle.

It is in fact a little tricky to construct a continuous family of operators from a continuous family of principal symbols but in fact we can prove homotopy type invariance without our doing so.

To see this note that if the pseudo-differential operators $P, Q$ have principal symbols such that

$$(1 - t)\sigma_m(P) + t\sigma_m(Q)$$

is elliptic for $0 \le t \le 1$ then the family

$$(1 - t)P + tQ$$

has the right principal symbols and thus $P, Q$ have the same index. So for a general family, using the fact that $S^*(X)$ is compact, we divide our homotopy into little pieces such that on each of these the straight homotopy works and homotopy invariance follows.

In fact, index is a pretty boring object for scalar differential operators. On manifolds of dimension bigger than two all operators have index zero. We won't prove this as its topology not analysis - it follows from the fact that all maps from $S^*(X)$ to $\mathbb{C} - \{0\}$ are null homotopic over coordinate patches and some use of Cech cohomology.

A consequence of this is that if an operator has trivial kernel then it is bijective which is certainly not obvious.

If one considers systems of operators or more generally operators on vector bundles then index is not necesarily zero.

## 12. Systems of Operators and Vector Bundles

We can consider a system of pseudo-differential operators on a manifold $X$ to be a matrix of scalar pseudo-differential operators.

**Definition 12.1.** *A pseudo-differential operator of order $m$ from*

$$C^\infty(X; \mathbb{C}^k) \to C^\infty(X; \mathbb{C}^l)$$

*is a map of the form*

$$(Pv)_i = \sum P_{ij} v_j$$

*$P_{ij}$ a scalar pseudo-differential operator of order $m$.*

*The principal symbol of $P$ is then just the the matrix of principal symbols, $\sigma_m(P_{ij})$.*

All the usual theorems go over and the matrix product of the principal symbols is principal symbol of the product. We thus define

**Definition 12.2.** *A system of pseudo-differential operators is elliptic if the principal symbol is invertible everywhere.*

And thus have

**Proposition 12.3.** *An elliptic system of pseudo-differential operators has a two-sided parametrix.*

*Proof.* Standard proof take $\sigma_{-m}(Q) = \sigma_m(P)^{-1}$ and then proceed as before. □

Clearly an analagous theorem will hold if the principal symbol is left (or right) invertible.

What we really want to do is consider operators between vector bundles which is not a large step from systems. A vector bundle a natural generalization of the situation where we have

$$\pi_0 : X \times \mathbb{R}^k \to X$$

with $X$ a smooth manifold and $\pi_0$ the projection onto $X$.

Recall that a vector bundle is a triple $(E, X, \pi)$ with $E$ a manifold of dimension $n + k$, $X$ a manifold of dimension $n$ and $\pi$ a projection such that the projection is locally trivial. That is $X$ can be covered by open sets $U_\alpha$ and there are diffeomorphisms

$$f_\alpha : \pi^{-1}(U_\alpha) \to U_\alpha \times \mathbb{C}^k$$

with $\pi f_\alpha = \pi_0$ and such that where two maps are defined, we have

$$f_\beta f_\alpha^{-1}(p, v) = (p, G_{\alpha\beta}(p)v)$$

with $G_{\alpha\beta}$ a smooth family of invertible matrices. This means that there is a natural linear structure on the fibres $E_p = \pi^{-1}(p)$ and so they are $k$ dimensional vector spaces.

One can then all the things to vector bundles one can do to vector spaces - e.g. tensor product $V \otimes W$, direct sum $V \oplus W$, homomorphisms $\mathrm{Hom}(V, W)$. One can also pull-back bundles by smooth maps. For a discussion of these operations see Atiyah [2]. By a section of a vector bundle $E$ over $X$, we mean a smooth map $f : X \to E$ such that $\pi f = \mathrm{Id}$. That is to each point, we assign a point in the fibre above it in a smoothly varying way, this class is denoted $C^\infty(X, E)$. We can define everything analagously for complex vector bundles.

We can define pseudo-differential operators on sections of vector-bundles.

**Definition 12.4.** *Let $E$ and $F$ be complex $C^\infty$ vector bundles over the manifold $X$ then a (classical) pseudo-differential operator of order $m$ from sections of $E$ to sections of $F$ is a continuous linear map*

$$A : C_0^\infty(X, E) \to C^\infty(X, F)$$

*such that for every open $U \subset X$ where $E, F$ are simultaneously trivialized then with respect to the trivializations, $A$ is given by a (classical) pseudo-differential system of order $m$. The class will be denoted $\Psi DO^m(X; E, F)$.*

It is of course enough to do this for a covering whilst assuming that operators between trivializations with disjoint bases are smoothing. A principal symbol is then defined before but will be a section of the pull-back of $\mathrm{Hom}(E, F)$ to $T^*(X)$. Locally, this just means it is a smooth matrix varying from point to point of $T^*(X)$ and will be homogeneous of order $m$ in $\xi$.

One can now repeat the entire theory of pseudo-differential operators for operators on sections of vector bundles and define Sobolev spaces of sections and prove all the same results. We leave the details to the enthusiastic reader. The important fact for us will be that an elliptic operator - ie an operator with invertible principal symbol - will be Fredholm as above. The index will then be an invariant of the homotopy class of the map as a function from

$$S^*(X) \to \mathrm{GL}(V, W)$$

where $\mathrm{GL}(V, W)$ is the bundle of invertible homomorphisms as a subset of $\mathrm{Hom}(V, W)$ above. Hence index theory!

## 13. The Hodge Theorem

So let's apply all this theory to do some differential geometry. Recall that over any smooth manifold, $X$, there are naturally the bundles of differential forms, $\Omega^k(X)$ and the natural operator $d$ between them.

Over any coordinate patch these bundles have a natural trivialization coming from the bases $dx_\alpha = dx_{\alpha_1} \wedge dx_{\alpha_2} \wedge \cdots \wedge dx_{\alpha_n}$ with $\alpha = (\alpha_1, \ldots \alpha_k)$ with $\alpha_j$ strictly increasing. That is locally a differential $k$-form can be written as a sum

$$\omega = \sum f_\alpha(x) dx_\alpha$$

with $\alpha$ of this form. And we have that

$$d\omega = \sum_{\alpha, j} \frac{\partial f_\alpha}{\partial x_j} dx_j \wedge dx_\alpha.$$

As $\wedge$ is anti-symmetric on one-forms it follows that $d^2 = 0$. What is the symbol of $d$? Well, the term $f_\alpha dx_\alpha$ is mapped to

$$\sum_j \frac{\partial f_\alpha}{\partial x_j} dx_j \wedge dx_\alpha$$

which has symbol given by mapping $dx_\alpha \to i \sum \xi_j dx_j \wedge dx_\alpha$. But a point in the cotangent bundle is a one-form that is a sum

$$\sum \xi_j dx_j$$

so the symbol of $d$ at $(x, \xi)$ is just

$$\mu \mapsto i\xi \wedge \mu.$$

The interesting thing about this symbol is that the kernel of the symbol of $d_k : \Omega^k \to \Omega^{k+1}$ is precisely the image of the previous one. To see this, observe that since everything is defined invariantly, we can pick coordinates such that the point $(x, \xi)$ is such that $\xi = (\xi_1, 0, \ldots, 0)$ saying something is in the kernel of the symbol is then just saying that the coefficient of $dx_\alpha$ is non-zero if and only if $\alpha_1 = 1$ but this is precisely the same condition as being in the image.

With this in mind, we define

**Definition 13.1.** *An elliptic complex (of order $m$) of vector bundles is a sequence $E_0, \ldots, E_k$ over a compact manifold $X$ together with pseudo-differential operators $P_j \in \Psi^m(X; E_j, E_{j+1})$ such that $P_j P_{j-1} = 0$ and the image of the symbol of $P_j$ is the kernel of the symbol of $P_{j+1}$ each $j$ with the first one injective and the last one surjective.*

   *(ie the symbols form an exact sequence)*

Note that if $k = 1$ then this says that the symbol is bijective everywhere that is the single operator is elliptic.

Provided each bundle is equipped with an inner product (which is always possible), we can use this structure to give an elliptic operator on each bundle. The inner product allows to define an adjoint operator $P_j^* \in \Psi DO^m(X : E_{j+1}, E_j)$. Putting these together we obtain operators on sections of $E_j$. We define

$$\Delta_j = P_j^* P_j + P_{j-1} P_{j-1}^* \in \Psi DO^{2m}(X; E_j, E_j)$$

for $j = 0, \ldots, k$. These are the desired elliptic operators. If $p_j$ is the symbol of $P_j$ then the symbol of $\Delta_j$ is $d_j = p_j^* p_j + p_{j-1} p_{j-1}^*$. If $d_j v = 0$ then

$$\langle d_j v, v \rangle = 0$$

and thus that

$$\langle p_j v, p_j v \rangle + \langle p_{j-1}^* v, p_{j-1}^* v \rangle = 0.$$

This means that $p_j v = 0$ which implies (by hypothesis) that $v = p_{j-1} w$ and

$$0 = \langle p_{j-1}^* v, w \rangle = \langle v, p_{j-1} w \rangle = \langle v, v \rangle.$$

So $v = 0$ and the ellipticity follows. The operators $\Delta_j$ then have elliptic, two-sided parametrices $T_j$ and are therefore Fredholm.

We can also use these parametrices to construct a *parametrix for the complex* that is to construct operators $Q_j$ such that

$$P_{j-1} Q_{j-1} + Q_j P_j - I$$

is smoothing for each $j$. To see this observe that

$$P_j \Delta_j = \Delta_{j+1} P_j$$

and therefore that

$$P_j T_j = T_{j+1} P_j$$

up to smoothing. Let $Q_j = T_j P_j^*$ and we have up to smoothing that

$$\begin{aligned}
P_{j-1} Q_{j-1} + Q_j P_j &= P_{j-1} T_{j-1} P_{j-1}^* + T_j P_j^* P_j \\
&= T_j (P_{j-1} P_{j-1}^* + P_j^* P_j) \\
&= I
\end{aligned}$$

The Generalized Hodge Decomposition Theorem is now very easy:

**Theorem 13.2.** *Suppose $P_j, E_j, j = 1, \ldots, k$ define an elliptic complex on a compact manifold $X$ then for all $j$*

$$C^\infty(X, E_j) = \ker \Delta_j \oplus P_{j-1} C^\infty(E_{j-1}) \oplus P_j^* C^\infty(E_{j+1})$$

*and*

$$\ker \Delta_j = \ker P_j \cap \ker P_{j-1}^*.$$

*Proof.* To prove the last statement, we have

$$0 = \langle \Delta_j v, v \rangle = \langle P_j v, P_j v \rangle + \langle P_{j-1}^* v, P_{j-1}^* v \rangle.$$

We know $\Delta_j$ is an elliptic self-adjoint operator so we have that

$$C^\infty(M) = \ker \Delta_j \oplus \Delta_j(C^\infty(M)).$$

But $\Delta_j = P_j^* P_j + P_{j-1} P_{j-1}^*$ so the decomposition follows.

We need to check the three subspaces are mutually orthogonal to complete the theorem. The orthogonality of the last two comes from $P_j P_{j-1} = 0$ as

$$\langle P_{j-1} v, P_j^* w \rangle = \langle P_j P_{j-1} v, w \rangle = 0.$$

The orthogonality of the image of $P_{j-1}$ with $\ker \Delta_j$ is just

$$\langle P_{j-1} v, w \rangle = \langle v, P_{j-1}^* w \rangle = 0$$

for

$$w \in \ker \Delta_j \subset \ker P_{j-1}^*$$

and the remaining relation is similar.                              $\square$

We can apply this to the de Rham complex to deduce things about the cohomology of a compact manifold. Recall that a $j$-form is closed if and only it is in the kernel of

$$d_j : C^\infty(X, \Omega^j(X)) \to C^\infty(X, \Omega^{j+1}(X)).$$

A $j$−form is exact if it is in the image of $d_{j-1}$. We make similar definitions for co-closed and co-exact in terms of the adjoint. The fact that $d^2 = 0$ means that the exact forms are closed and the co-exact are co-closed. A form that is closed and co-closed is said to be harmonic. The theorem above says that any form is the sum of a harmonic form, an exact form and a co-exact form.

Recall that the $j$-th de Rham cohomology of a manifold is defined to be the space of closed $j$-forms quotiented by the space of exact $j$-forms. We claim that this is naturally isomorphic to the space of harmonic $j$-forms. This will follow trivially from Theorem 13.2 if we show that the space of closed forms is equal to the space of exact forms plus the space of harmonic forms. Now a form is harmonic if and only if it is both closed and co-closed - same proof as ellipticity of $\Delta_j$. Using the decomposition, it is thus enough to show that any non-zero co-exact form is never closed, but

$$\langle dd^*v, v \rangle = \langle d^*v, d^*v \rangle$$

which can only be zero if $d^*v = 0$ and the result follows. As $\Delta_j$ is Fredholm, we have that the space of harmonic forms is finite-dimensional so the de Rham cohomology of a compact manifold is finite dimensional - a non-obvious fact!

NB there is a slight flaw in the above in that it was done for complex rather than real forms but we can commute taking real parts through the final results to prove that they hold in the real case too.

We proved all this while using any inner product on each form bundle - if one is more careful and defines all the inner products using the Riemannian structure of the tangent bundle then one can also use this to prove Poincare duality. To do this we have to use several natural pairings which are at our disposal.

First there is a natural pairing between the exterior algebra of a vector space and the exterior algebra of its dual.

$$\Lambda^j(V) \times \Lambda^j(V^*) \to \mathbb{R}$$

Defined on basis elements by $\alpha = \mu_1 \wedge \cdots \wedge \mu_k$, $\beta^* = f_1 \wedge \cdots \wedge f_k$ by

$$(\alpha, \beta^*) = \det(f_i(\mu_j)). \tag{9}$$

This can be extended by linearity and by considering the pairing the wedge products of an element of a basis and the wedge products of the associated dual basis, it is non-singular. It is also intrinsic. So we have that

$$\Lambda^j(V^*) \sim (\Lambda^j(V))^*.$$

If we pick a volume element, that is a non-zero element $\omega$ of $\Lambda^n(V)$, then we have a pairing

$$\Lambda^j(V) \times \Lambda^{n-j}(V) \to \mathbb{R}$$

$$(\alpha, \beta) \mapsto c(\alpha, \beta)$$

where

$$\alpha \wedge \beta = c(\alpha, \beta)\omega.$$

This is non-singular also. So given a volume form we have

$$\Lambda^j(V) \sim (\Lambda^{n-j}(V))^*. \tag{10}$$

So combining (9) and (10), we have

$$\Lambda^j(V) \sim \Lambda^{n-j}(V^*). \tag{11}$$

However, we want an isomorphism with $\Lambda^{n-j}(V)$ rather than $\Lambda^{n-j}(V^*)$. To do this, we assume we have an inner product, $B$, on $V$ and thus a natural isomorphism, $L_B$, from $V$ to $V^*$. Because the exterior algebra is a natural object this means that we have an isomorphism

$$L_B^{-1} : \Lambda^{n-j}(V^*) \sim \Lambda^{n-j}(V).$$

Note also that once we have chosen an orientation of $V$, the bilinear form can be used to fix the choice of volume form - an orientation is a choice of volume form up to a positive multiple. We define an *oriented basis* of $V$ to be a basis, $e_l$, such that $e_1 \wedge \cdots \wedge e_n$ is in the orientation class. Any orthonormal oriented basis will then give the same volume form.

So having chosen an inner product and an orientation, we define

$$* : \Lambda^j(V) \to \Lambda^{n-j}(V)$$

to be the given isomorphism from $\Lambda^j(V)$ to $\Lambda^{n-j}(V)$ followed by $L_B^{-1}$. Tracing this through on any oriented orthonormal basis, $v_1, \ldots, v_n$, we have

$$*(v_1 \wedge \ldots v_j) = v_{j+1} \wedge \cdots \wedge v_n, \tag{12}$$

which allows us to compute $*$ on any element and we have that

$$** = (-1)^{j(n-j)} \, \mathrm{Id}$$

on elements of order $j$.

We can now define ( yet another ) pairing,

$$(\alpha, \beta) = \alpha \wedge *\beta = c(\alpha, \beta)\omega_B \tag{13}$$

where $\omega_B$ is the volume form coming from the inner product.

We now want to do all this globally. Suppose we have an oriented, compact, Riemannian manifold, $X$. Then the above construction will work with $V = T_p^*(M)$ for any $p \in M$. So we have an inner product on each form bundle. What is the adjoint of $d$, $\delta$, with respect to this inner product on the $j$-th form bundle?

$$\delta = d^* = (-1)^{j+1+j(n-j)} * d * .$$

To see this observe that by Stokes' theorem

$$\int d(\alpha \wedge *\beta) = 0$$

and compute $d(\alpha \wedge *\beta)$.

What does all this do for us? Well since

$$\Delta = d\delta + \delta d$$

it commutes with $*$ up to sign, which means that $*$ is a bijection from the space of harmonic forms in dimension $j$ to that in dimension $n-j$, ie Poincare duality:

$$H^{n-j}(X) = H^j(X).$$

It is also worth mentioning that the isomorphism of de Rham cohomology with the space of harmonic forms means that by making assumptions on the metric one can prove things about the cohomology.

## 14. EXERCISES

*Question* 1. Construct a function, $\phi$, with the following properties:

1. $\phi$ is a smooth, real-valued function on $\mathbb{R}^n$,
2. $\phi = 1$ in $|x| \leq 1$,
3. $\phi = 0$ in $|x| \geq 2$,
4. $\phi \geq 0$.

*Use the function*

$$f(x) = \begin{cases} 0 & x \leq 0, \\ e^{-\frac{1}{x}} & x > 0. \end{cases}$$

*Question* 2. Show $C_0^\infty(\mathbb{R}^n)$ is dense in $\mathcal{S}(\mathbb{R}^n)$ in the sense that given $\phi \in \mathcal{S}(\mathbb{R}^n)$ there exists a sequence $\phi_n \in C_0^\infty(\mathbb{R}^n)$ such that

$$||\phi - \phi_n||_{\alpha,\beta} \to 0, \ \forall \alpha, \beta.$$

Deduce that if $u \in \mathcal{S}'(\mathbb{R}^n)$ is zero on $C_0^\infty(\mathbb{R}^n)$ then it is zero on $\mathcal{S}(\mathbb{R}^n)$. Observe that $\mathcal{S}'(\mathbb{R}^n)$ then naturally forms a subspace of $\mathcal{D}'(\mathbb{R}^n)$.

*Question* 3. (Partitions of Unity) Suppose $U_\alpha$ is a collection of open subsets of $\mathbb{R}^n$ such that $\cup U_\alpha = \mathbb{R}^n$ and $\bar{U}_\alpha$ is compact for each $\alpha$ and any infinite intersection of $U_\alpha$ is empty. Show that there exists a collection of smooth functions $\phi_\alpha$ such that $\text{supp}(\phi_\alpha) \subset U_\alpha$ and

$$\sum_\alpha \phi_\alpha = 1.$$

*If you know about paracompactness, what can one do when we do not assume the intersections are locally finite.*

*Question* 4. Suppose $u \in \mathcal{D}'(\mathbb{R}^n)$ and $\text{supp}(u)$ is compact. Show that the pairing

$$< u, e^{-ix.\xi} >$$

makes sense for any $\xi \in \mathbb{C}^n$ and that it results in a function which is smooth in $\xi$. Show also that it satisifies the Cauchy-Riemann equations in each $\xi_j$.

*Question* 5. Euler's Relation

Suppose $f$ is a smooth function on $\mathbb{R}^n - \{0\}$ and that it is (positively) homogeneous of order $m$. ie

$$f(\lambda x) = \lambda^m f(x)$$

for all $\lambda > 0$. Show that

$$\left( \frac{\partial}{\partial x} x - m \right) f(x) = 0.$$

Prove the converse.

Show that there is a one-one correspondence between the set of homogeneous functions of order $m$ and the smooth functions on $S^{n-1}$.

Suppose $f$ is smooth on $\mathbb{R}^n$, show that $f$ is a homogeneous polynomial.

*Question* 6. Suppose $f(x, \xi)$ is a smooth function on $\mathbb{R}_x^n \times (\mathbb{R}_\xi^k - \{0\})$ which is homogeneous of degree $\mu \in \mathbb{C}$ in $\xi$ and $\phi$ is as in $Q1$. Show that

$$f(x, \xi)(1 - \phi)(\xi)$$

is a symbol of order $\Re \mu$.

*Question* 7. The Borel Lemma

The point of this question is to show that given any sequence $\{c_n\}$ there exists a function whose Taylor series is $\{c_n\}$. i.e. there exists a smooth function $f \in C^\infty(\mathbb{R})$ such that

$$\partial_x^j f(0) = j! c_j.$$

Let $\phi$ be as in Q1.

1. Show that if $t_j \to \infty$ then the sum

$$\sum_{j=0}^\infty \phi(t_j x) c_j x^j$$

   converges for each $x$.
2. Show that if $\{t_j\}$ is picked correctly then the sum will converge uniformly to a continuous function.
3. Show there exists $\{t_j^k\}$ such that the sum of the $k^{th}$ derivatives converges uniformly .
4. Use a diagonalization argument to show that the same $\{t_j\}$ can be used for all $k$ and deduce that then

$$\sum c_j \phi(t_j x) x^j$$

   is the desired function.
5. Need the sum $\sum_{j=0}^\infty c_j x^j$ exist for any non-zero $x$ ?

*Question* 8. Asymptotic Summation Show that if $a_j \in S^{m-j}(\mathbb{R}^n; \mathbb{R}^k)$ then there exists $a \in S^m(\mathbb{R}^n; \mathbb{R}^k)$ such that

$$a - \sum_{j=0}^N a_j \in S^{m-N} \text{ for all } N.$$

  *Repeat the proof of the Borel Lemma but do it about infinity instead of about 0.*

*Question* 9. Find a formula for the total symbol of the composite of two differential operators in terms of their total symbols. Identify a formula for the term of each homogeneity. (*reduce to simple cases*)

*Question* 10. Define a distribution, $u$ on $\mathbb{R}^n$, to be homogeneous of order $m$ if

$$\left( x \frac{\partial}{\partial x} - m \right) u = 0.$$

Observe that a homogeneous distribution is tempered and show that its Fourier transform is also homogeneous.

NB $x \frac{\partial}{\partial x} = \sum\limits_{j=1}^{n} x_j \frac{\partial}{\partial x_j}$

For the brave - when can a homogeneous smooth function on $\mathbb{R}^n - 0$ be extended to a homogeneous distribution on $\mathbb{R}^n$?

*Question* 11. Suppose $a \in S^m(\mathbb{R}_x^n; \mathbb{R}_\theta^k)$. Show that $a \in S_{cl}^m$ then

$$\prod_{j=0}^{N-1} \left( \theta \frac{\partial}{\partial \theta} - (m - j) \right) a(x, \theta) \in S^{m-N} \forall N.$$

*Question* 12. Show that if $P$ is a pseudo-differential operator and $Q$ is smoothing and at least one is proper then $PQ, QP$ is smoothing.

*Question* 13. If $P$ is a properly supported operator with Schwartz kernel $K$ and defines maps from $\mathcal{D}'(\mathbb{R}^n) \to \mathcal{D}'(\mathbb{R}^n)$ and from $C^\infty(\mathbb{R}^n) \to C^\infty(\mathbb{R}^n)$, show that

$$\operatorname{singsupp}(Pu) \subset \operatorname{singsupp}(K) \circ \operatorname{singsupp}(u).$$

*Decompose $P$ and $u$ into pieces supported near and far from their singular supports and then use the relation for ordinary supports.*

*Question* 14. If $P$ is a classical pseudo-differential operator with positive principal symbol, use a symbolic argument to show there exists $Q$ a classical, pseudo-differential operator such that

$$P - Q^2 \in \Psi DO^{-\infty}.$$

*Question* 15. Show that if $a \in S^m$ is classical then its asymptotic expansion in homogeneous terms is unique.

*Question* 16. The Weyl Quantization

Show that if $P \in \Psi DO_{cl}^m(\mathbb{R}^n)$ then $P$ can be written in the form

$$\int e^{i<x-y,\xi>} a \left( \frac{x+y}{2}, \xi \right) d\xi$$

up to smooth terms for some classical symbol $a$ and express $a$ in terms of the left symbol.

*Question* 17. Exhibit a distribution, $u$, satisfying

$$\left(D_t^2 - D_x^2\right) u(x,t) = 0$$

which is not smooth.

*Question* 18. Classify the solutions of $D_{x_1} u = 0$ on $\mathbb{R}^n$ and classify the possible singular supports.

*Question* 19. Give an example of a Schwartz kernel which maps $C_0^\infty$ into $\mathcal{D}'$ but not $C^\infty$.

*Question* 20. If $P$ is a pseudo-differential operator on $\mathbb{R}^n$ then

$$P - P^*$$

is smoothing if and only if the Weyl symbol of $P$ is real. (see last problem sheet.)

*Question* 21. In what Sobolev classes are the following? give the answer for both the $H^s$ and $H_{loc}^s$

1. $\delta(x)$ on $\mathbb{R}^n$
2. $H(x) = \begin{cases} 1 & x \leq 0 \\ 0 & x < 0 \end{cases}$
3. $1/x$ on $\mathbb{R}$ (defined by principal values)
4. the characteristic function of the unit disc in $\mathbb{R}^2$.
5. The distribution obtained by integrating a function over the unit sphere.

*Question* 22. Compute the wavefront set of each of the distributions in Q21//. Let $g^{ij}(x)$ be a positive definite matrix smoothly varying with x and define

$$\Box = D_t^2 - \sum_{i,j} g^{ij}(x) D_{x_i} D_{x_j}.$$

Suppose $\Box u$ is smooth, give an upper bound on $\mathrm{WF}(u)$.

*Question* 23. If $P_j \in \Psi DO_{cl}^{m_j}$ show that

$$[P,Q] = PQ - QP \in \Psi DO_{cl}^{m_1+m_2-1}$$

and compute the principal symbol. Show that principal symbols form a Lie Algebra under this operation. (if you know what a Lie algebra is.)

*Question* 24. Try to prove the coordinate invariance of $H_c^s$ directly.

If $p_m(x, \xi)$ is a homogeneous function on $T^*(\mathbb{R}^n) - \{0\}$, define the Hamiltonian of $p_m$ to be the vector field

$$\sum \frac{\partial p_m}{\partial \xi_j} \frac{\partial}{\partial x_j} - \frac{\partial p_m}{\partial x_j} \frac{\partial}{\partial \xi_j}.$$

Show that the Hamiltonian is coordinate invariant. Show also that the integral curves of the Hamiltonian preserve the set $p_m = 0$. Show that multiplying $p_m$ by a non-zero function does not change the integral curves within $p_m = 0$.

*Question* 25. Construct a distribution on $\mathbb{R}^n$ of which the wavefront set is a single ray.

*Question* 26. What is the wavefront set of $H(x_1)\delta(x'')$ on $\mathbb{R}^n$?

*Question* 27. Prove the continuity of zeroth order non-classical pseudo-differential operators on $L_c^2$.

*Question* 28. If $P \in \Psi DO_{cl}^m$ with positive principal symbol prove that $P$ has an $m^{th}$ root up to smoothing operators.

*Question* 29. Show that the intersection of the characteristic varieties of vector fields tangent to a submanifold is equal to the submanifold's conormal bundle. Deduce an upper bound for the wavefront set of a pseudo-differential operator.

*Question* 30. Suppose $P$ is a pseudo-differential operator on $\mathbb{R}^n$, define

$$P' : C_c^\infty(\mathbb{R}^n \times \mathbb{R}) \to C^\infty(\mathbb{R}^n \times \mathbb{R})$$

by treating the extra variable as a parameter. Show $P'$ is a pseudo-differential operator if and only if $P$ is a differential operator.

*Question* 31. We shall say $u \in \mathcal{D}'(X)$ is in $H_{loc}^s$ at $(x_0, \xi_0)$ if there exists a zeroth order pseudo-differential operator micro-elliptic at $(x_0, \xi_0)$ such that

$$Pu \in H_{loc}^s(X).$$

Show that if $u \in H_{loc}^s$ at every point $(x_0, \xi_0)$ then $u \in H_{loc}^s$.

Show also that if $Q$ is order $m$ and $u \in H_{loc}^s$ at $(x_0, \xi_0)$ then $Qu \in H_{loc}^{s-m}$ at $(x_0, \xi_0)$. Deduce that if $Q$ is micro-elliptic at $(x_0, \xi_0)$ and $Qu \in H_{loc}^s$ at $(x_0, \xi_0)$ then $u \in H_{loc}^{s+m}$ at $(x_0, \xi_0)$.

Define $s_u(x, \xi)$ to the supremum of $s$ such that $u \in H^s_{loc}$ at $(x, \xi)$. Is $s_u$ continuous or semi-continuous?

Does $s_u(x, \xi) = +\infty$ imply that $(x, \xi) \notin \mathrm{WF}(u)$ ?

*Question* 32. Show that if $V$ is a vector bundle then there exists a smooth metric on $V$ ie an inner product on each fibre which varies smoothly from point to point.

*Question* 33. Show that there is natural pairing

$$\Gamma(|\Lambda|^s, X) \times \Gamma(|\Lambda|^{1-s}, X) \to \mathbb{C}$$

given by integration.

*Question* 34. If $\Lambda^n$ is the line bundle of n-forms on a manifolds, show there is a natural trivilization over each coordinate patch and compute the transition functions.

*Question* 35. Let $\square = D_t^2 - \Delta$ be the flat wave equation on $\mathbb{R} \times \mathbb{R}^2$. Let $P$ be the parabola $y = x^2$ and let $u$ be a distribution such that

$$\square u = 0.$$

Suppose

$$\mathrm{WF}(u) \cap \{t = 0\} \subset \{(x, y, \xi, \eta, |\xi, \eta|) : (x, y, \xi, \eta) \in N^*(P)\}.$$

Compute the wavefront set of $u$ and sketch the projection of $\mathrm{WF}(u)$ onto $(x, y)$ for $t$ at important times.

*Question* 36. Show that if $\square u$ is smooth and singsupp$(u)$ is compactly supported then $u$ is smooth.

*Question* 37. Give an example of a differential operator which has periodic bicharacteristics.

*Question* 38. * Suppose $P$ is a zeroth order elliptic, pseudo-differential operator on the circle with principal symbol equal to

$$\sigma_0(P) = \begin{cases} a_+(x) & \xi > 0 \\ a_-(x) & \xi < 0. \end{cases}$$

Compute the index of $P$ in terms of the winding numbers of $a_\pm$ as maps from $S^1$ to $\mathbb{C} - \{0\}$.

Deduce that the result holds for $P$ of any order.

*Question* 39. If $\omega \in S^{n-1}$ show that the distribution

$$\delta(t - x.\omega)$$

is a solution of

$$\Box u = 0$$

and compute its wavefront set.

*Question* 40. Are the left and right shift operators on $l_2(\mathbb{N})$ Fredholm? If so what are their indexes?

*Question* 41. Let $\mathcal{M}$ be a subset of $\Psi DO_{cl}^1$ and define

$$I^{(s)}(\mathcal{M}) = \{u : P_1 \ldots P_k u \in H_{loc}^s, \forall k, P_j \in \mathcal{M}\}.$$

Show that if $u \in I^{(s)}(\mathcal{M})$ then

$$\mathrm{WF}(u) \subset \bigcap_{P \in \mathcal{M}} \sigma_1(P)^{-1}(0).$$

Show that if $\overline{\mathcal{M}}$ is the smallest subset of $\Psi DO_{cl}^1$ containing $\mathcal{M}$ which is closed under summation, multiplication by zeroth order operators and the taking of commutators then

$$I^{(s)}(\overline{\mathcal{M}}) = I^{(s)}(\mathcal{M}).$$

*Question* 42. Let $S_{\rho,\delta}^m$ be the class of symbols $a$ such that

$$|D_x^\beta D_\xi^\alpha a(x,\xi)| \leq C_{\alpha,\beta,K} <\xi>^{m-|\alpha|\rho+|\beta|\delta}$$

for $x$ in $K$ compact with $\rho,\delta$ between 0 and 1. For which $\rho,\delta$ does the theory of pseudo-differential operators on $\mathbb{R}^n$ work? For which $\rho,\delta$ does the theory work on manifolds?

*Question* 43. Let $X$ be a compact manifold and suppose $\{\Gamma_j\}$ is a collection of open cones in $T^*(X) - 0$ which form an open cover. Show there exists a collection of zeroth order pseudo-differential, operators, $\{P_j\}$, such that

$$\sum_j P_j = \mathrm{Id}$$

and

$$\mathrm{WF}'(P_j) \subset \Gamma_j.$$

Find out what a sheaf is and deduce that the space of distributions on $X$ modulo the smooth functions form a sheaf over the cosphere bundle.

*Question* 44. Show that if $a_j \in S^{m-j}$ and $a \in S^{m'}$ are such that there exists a sequence $N_j$ such that $N_j \to -\infty$ and

$$a - \sum_{j<N} a_j \in S^{N_j}$$

then $a \in S^m$ and

$$a - \sum_{j<N} a_j \in S^{m-j}$$